\documentclass[10pt,a4paper]{article}

\usepackage{amsmath,amssymb,latexsym,color,dsfont,mleftright}
 \usepackage{authblk,mathtools, nccmath,subcaption}
 \usepackage{graphicx,float}
 \usepackage[square,numbers]{natbib}
 \usepackage{tikz-cd}
  \usepackage{bbm}
 \usepackage[cal=boondox]{mathalfa}

 \usepackage{tikz,xcolor,hyperref}
   \usepackage{cleveref}
 \usepackage{orcidlink}

\newtheorem{theorem}{Theorem}
 
\newtheorem{conjecture}{Conjecture}
\newtheorem{lemma}{Lemma}
\newtheorem{proposition}{Proposition}
\newtheorem{definition}{Definition}

\newtheorem{problem}{Problem}
\newtheorem{remark}{Remark}

\newcommand{\lqqd}{\hfill{$\Box$}\bigskip}
\newcommand{\CC}{\mathbb{C}}

\newcommand{\NN}{\mathbb{N}}

\newcommand{\RR}{\mathbb{R}}

\newcommand{\LL}{\mathcal{L}}

\newcommand{\dgr}[1]{\mbox{{ \textrm{deg}\/}}[#1]}
\newcommand{\supp}[1]{\mbox{{\textrm{supp}\/}}[#1]}

\newcommand{\dsty}{\displaystyle}
\newcommand{\proof}{\bigskip\noindent {\sc Proof.  }}

\def\bbbuildrel#1_#2{\mathrel{
\mathop{\kern 0pt#1}\limits_{#2}}}

\title{An asymptotic expansion of  eigenpolynomials for a class of linear differential operators.}

\author{Jorge A. Borrego-Morell   \\
Universidade Federal do Rio de Janeiro\\
Campus UFRJ Duque de Caxias, Brazil}
 
\begin{document}

%%% ----------------------------------------------------------------------
\maketitle
%%% ----------------------------------------------------------------------
 
\begin{abstract}
Consider an   $M$-th order linear  differential operator, $M\geq 2$, 

$$
\LL^{(M)}=\sum_{k=0}^{M}\rho_{k}(z)\frac{d^k}{dz^k},
$$
where $\rho_M $ is a monic complex polynomial such that $\dgr{\rho_M }=M$ and  $(\rho_k)_{k=0}^{M-1}$ are complex polynomials such that $\dgr{\rho_k }\leq k, 0\leq k \leq M-1$. It is known that the zero counting measure  of its eigenpolynomials  converges in the weak star sense  to a measure $\mu$. 
We obtain an asymptotic expansion of  the  eigenpolynomials of $\LL^{(M)}$ in compact subsets out the support of $\mu$. In particular, we solve  a conjecture posed in \cite{masshap01}.

\end{abstract}

%\noindent {\it MSC:} 33C45, 33C47, 42C05, 34A05.\\


\section{Introduction}
Let $M \geq 2$ be an integer and consider the linear ordinary differential operator of order $M$
\begin{equation}\label{DiffOper}
\LL^{(M)}=\sum_{k=0}^{M} \rho_k(z) \frac{d^k}{dz^k},
\end{equation}
where  $\dsty \rho_k(z)=\sum_{j=0}^k\rho_{k,j}z^j,k=0,1, \dots, M$  are  polynomials in $z$ such that $\dgr{\rho_k } \leq k$ for $k=0,1, \dots, M$ and at least one of them, say $\rho_{k^*} $, is exactly of degree $k^*$,  $\rho_M  $ being  a monic polynomial. Note that if $f_n$ is a polynomial of degree $n$ then $\dgr{\LL^{(M)}[f_n]}\leq n$, linear differential operators that satisfy this property of invariance are called \emph{exactly-solvable} (cf. \cite{Turb92b}) and can be split in two families: non-degenerate, if the leading term of the operator satisfies $\dgr{\rho_M } =M$ and degenerate, if $\dgr{\rho_M }<M$. In particular, Bochner--Krall operators, defined as linear differential operators having a sequence of  eigenpolynomials  which are  orthogonal polynomials,  belong to the class of exactly solvable linear differential  operators,  see \cite{Ev01,HST}  for further details.

 General properties  of the eigenpolynomials of   exactly solvable linear differential  operators have been previously studied in e.g. \cite{Sh33} and for operators of the form $\dsty \LL^{(M)}=\sum_{k=0}^{M} \rho_k(z) \frac{d^k}{dz^k}$  in \cite{masshap01}. Further studies of  exactly solvable operators in general,  mainly motivated by the Bochner problem,  can be found in  \cite{BerRull02,Berg07,BeRuSh04}.

For a polynomial  $P_n $   of degree $n$, denote by  $\dsty \mu_n=\frac{1}{n}\sum_{P_n(z_{k,n})=0}\delta_{z_{k,n}}$ the normalized counting measure of the zeros of  $P_n $. For a given sequence $(P_n)_{n=0}^{\infty}$ of polynomials, if the  limit $\dsty\mu = \lim_{n\to \infty}\mu_n$ exists  (in the sense of the weak convergence), then $\mu$ is called the asymptotic zero-counting measure of the sequence.

It is known that for the operator \eqref{DiffOper} and  $n$ large enough,  there exists a unique eigenpolynomial $Q_n$ of degree $n$ \cite[Th.1]{Sh33}, \cite[Th.1]{BerRull02} and its  eigenvalue is  given by,
\begin{eqnarray}\label{lambn} 
\lambda_n=\sum_{0\leq k\leq \min(M,n)}\rho_{k,k}(n)_k,
\end{eqnarray}
where $(n)_0=1, (n)_k=n(n-1)\cdots (n-k+1), k\in \NN$ is the falling factorial, cf. \cite[p.6]{comtet74}. Moreover, the zero counting measures   of the eigenpolynomials converge weakly to a measure $\mu$ with support, denoted through this paper as $\supp \mu $,  contained in the convex hull of the zeros of $\rho_{M} $. This support  is finite, connected, and  its complement is also connected, cf. \cite[Ths. 3,4]{BerRull02}.

 The study  of asymptotic properties of a system of polynomials is an important topic in the theory of orthogonal polynomials and approximation theory.  The objective  of this  manuscript is the proof of an asymptotic formula   for   the  monic polynomial eigenfunctions $Q_n$ of $\LL^{(M)}$ in compact subsets of $\CC\setminus \supp \mu$ in the sense of Poincar\'e.  The main tool is the  WKB approximation,  named   after the  physicists Wentzel \cite{wentzel}, Kramers \cite{kramers}, and Brillouin \cite{Brillouin}, who  efficiently  used this technique in the study of asymptotic solutions of linear second--order differential equations in  quantum physics.  The method was known before to Liouville \cite{Liov}, Green \cite{Green},  and other authors, see \cite[Ch. XIII]{Ding}, \cite[Ch.1]{wasow2012}, and \cite[p.228]{Olv74} for some historical remarks. 

 In general,  the  
infinite series that define WKB solutions do not converge in the
usual sense, see  \cite[p.365 ]{Olv74} for an example. The use of Borel resummed WKB solutions to overcome this situation was first suggested  in \cite{BendWu},  developed later   in \cite{Vor}, and further in a series of publications  \cite{AKT91,AKT09,BoSc02,DP99,kawai2005algebraic,N21}.  For  higher--order  linear  differential equations this method  has   also been studied considerably,  see   \cite{AokiKawTak94,Fed93,kawai2011exact,MT17,TAK17} and the references within. The basic idea behind is the analysis  of  the Stokes geometry which describes regions in the complex plane  (sometimes called Stokes regions) where WKB solutions are Borel summable and the connection formulas which describe relations between the Borel sums of WKB solutions in different Stokes regions. In \cite{MT17} the authors  study the  structure of the Stokes geometry of the WKB solutions of a linear $n$th-order differential operators  of the form 
$$ P\phi(z,\eta)=\sum_{j=0}^na_j(z)\left(\eta^{-1}\frac{\partial}{\partial z}\right)^j\phi(z,\eta)=0,$$
where $a_j$ are polynomials, $a_n$ a  non--zero complex constant, and $\eta>0$ is a real large parameter. By reducing    higher--order  linear ordinary differential equations to second--order   via  the  middle convolution operation introduced in \cite{K96} and developed in   \cite{DR07,O12},  the authors obtain some information about  the Stokes geometry    of the  WKB solutions of $P$. For higher--order equations (i.e., $n\geq 3 $),   the complete structure of the Stokes geometry is not fully understood yet. In this manuscript we will not use the approach of the Borel resummed solutions.  However, it would be of interest the application of this theory and  will be done in a further publication.

To obtain  the asymptotic expansion,   we   use instead   a result due to Sibuya \cite[Th.XII-1-2 p.374]{HSib99} and \cite{Si58} on the existence of asymptotic solutions in the sense of Poincar\'e of a perturbed system of  first--order linear differential equations.

\subsection{Main notation and statement of the results}

Let  ${\bf B}_{n,k}$ denote  the (exponential) partial Bell polynomials  \cite[(3a) (3d), p. 134]{comtet74}
\begin{eqnarray*}
\nonumber {\bf B}_{n,k}\left(x_1,\ldots,x_{n-k+1}\right)&=&
\sum   \frac{n!}{c_1!c_2!\cdots (1!)^{c_1}(2!)^{c_2}\cdots} x_1^{c_1}\cdots x_{n-k+1}^{c_{n-k+1}}, \\
\nonumber {\bf B}_{0,0}&=&1,
\end{eqnarray*}
where the summation is taking over all integers $c_1,c_2,\cdots \geq 0$ such that 
\begin{eqnarray*}
\nonumber c_1+2c_2+3c_3\cdots &=&n,\\
\nonumber c_1+c_2+c_3\cdots   &=&k.
\end{eqnarray*}

Let  ${\bf Y}_n={\bf Y}_n(x_1,\ldots,x_n), {\bf Y}_0=1$  be the (exponential) complete Bell polynomials  \cite[(3d), p. 134]{comtet74} defined by the generating function  in an infinite number of indeterminates  $x_1,x_2,\ldots $  given by
$$\exp \left(\dsty \sum_{m\geq 1}x_m\frac{t^m}{m!}\right)=1+\sum_{n\geq 1}{\bf Y}_n(x_1,\ldots,x_n)\frac{t^n}{n!}.$$

 Denote  by $\Delta(z^{\prime},r)$  the open disk $|z-z^{\prime}|<r$, and define also  the  open sector     $\mathcal{G}(r,\alpha)=\{\epsilon\in \CC: |\arg[\epsilon]|<\alpha, 0<|\epsilon|<r\}$. Given  an open set $U\subset\CC$, we denote by  $\mathcal{H}(U)$  the space of analytic functions in   $U$. Denote also the topological operators  $\mathring{U}, \partial U, \overline{U}$ and $U^c$  as the  interior, boundary, closure, and complement relative to $\CC$  respectively of $U$.   If  $\gamma$ is a closed Jordan curve, we denote by $int(\gamma)$ and $ext(\gamma)$  the  bounded and unbounded components respectively.   Let $\tau $ be any oriented Jordan arc, then   $\tau^+$ will  denote  the side on the left  of $\tau$. Define    $\Omega=\CC\setminus\supp \mu$, for a given Jordan arc    $\tau$   connecting  $\infty$ and  $\omega\in \supp \mu$ we will denote by $\Omega_{\tau}$ the open connected set $\Omega\setminus \tau$. Denote by $ w_1$ the branch of $\left(\sqrt[M]{\rho_M(z)}\right)^{-1}$ in $\Omega$ which tends to $\dsty \frac{1}{z}$ near $\infty$. By \cite[Th.3]{BerRull02} the set $\supp\mu$ is a finite  tree, notice that $w_1$ is then an analytic function out of $\supp\mu$.  Define also 
\begin{eqnarray}\label{wj}
w_j(z)=e^{\frac{2(j-1)\pi\imath}{M}}w_1(z), z\in \Omega, j=2,\ldots,M.
\end{eqnarray}
Let $\tau$ be any Jordan arc  from  $-\infty$ to $\omega\in \supp \mu$ such that $]-\infty,p]\subset \tau$, for some $p\in \RR$. Define  $\Phi_0 $ as  the  primitive of $w_1$ in $\Omega_{\tau}$   such that  $\dsty\lim_{  
 z\rightarrow \infty  } \Phi_0(z)-\ln z=0$ and    $\Phi_1 $ as  the primitive in $\Omega_{\tau}$ of the function
$$\dsty    \frac{(M-1)\rho_{M}^{\prime}(z)}{2M \rho_{M}(z)}-\frac{\rho_{M-1}(z)}{M\rho_{M}(z)}, $$
such that  $\dsty\lim_{
 z\rightarrow \infty  } \Phi_1(z)-\left(\frac{M-1}{2}-\frac{\rho_{M-1,M-1}}{M}\right)\ln z=0$.   The main result can be stated as follows,

\begin{theorem}\label{Main}
Let  $\LL^{(M)}, M\geq 2$ be a non--degenerate exactly solvable operator and $Q_{n} $ be  the $n$--th monic eigenpolynomial of $\LL^{(M)}$, 
then as $n\rightarrow\infty$
 $$\dsty Q_n(z)\sim e^{\left(\dsty n\Phi_{0}(z)-\left(\frac{M-1}{2}-\frac{\rho_{M-1,M-1}}{M}\right)\Phi_{0}(z)+\Phi_{1}(z)\right)}\left(1+\frac{C_1(z)}{n}+\frac{C_2(z)}{n^2}+\ldots\right),
$$
uniformly in compacts subsets  $K\subset \Omega$, here $C_j\in \mathcal{H}(\Omega), j\geq 1$ and the symbol $\sim$ has the same meaning as in \cite[(7.03) p.16]{Olv74} for an  asymptotic expansion in the sense of Poincar\'e.
\end{theorem}

\subsection{Overview of the proof and structure of the manuscript}

Consider the differential equation depending on the parameter $\epsilon$
\begin{eqnarray*}
\LL^{(M)}{[v]}(z,\epsilon)- \frac{v(z,\epsilon)}{\epsilon^{M}}=0, \quad \epsilon \in \CC \setminus \{0\},
\end{eqnarray*}
where $\LL^{(M)}$ is given by \eqref{DiffOper}. Writing 
\begin{eqnarray}\label{DiffEqn01}
v^{(M)}(z,\epsilon)+\sum_{k=1}^{M-1} \frac{\rho_k(z)}{\rho_M(z)} v^{(k)}(z,\epsilon)- \frac{v(z,\epsilon)}{\epsilon^{M} \rho_M(z) }=0,
\end{eqnarray}
where   $\dsty \rho_k(z)=\sum_{j=0}^k\rho_{k,j}z^j$, we see that   $Q_n $ is a polynomial solution of \eqref{DiffEqn01} with $\dsty \epsilon_n=\frac{1}{ \sqrt[M]{\lambda_{n}-\rho_{0,0}}}$. Reciprocally, if for some quantity $\epsilon_n$, $v_n$ is a polynomial solution of \eqref{DiffEqn01} of degree $n$, then $v_n$ is a polynomial eigenfunction of the exactly solvable differential operator \eqref{DiffOper} with eigenvalue $\lambda_{n}=\rho_{0,0}+\epsilon_n^{-M}$.

\begin{remark}
For brevity,   we will   use   the superscript $(k)$   for  the $k$th-order partial derivative of $v$ with respect to  $z$.  

\end{remark}

The main result of  Section \ref{Li} is Proposition \ref{GL}  which guarantees the  existence of the exponential series that define the WKB solution when $z$ varies in a disk and   $\epsilon$ varies in a sectorial region. More precisely,   there exist  $\alpha,r,\eta^{\prime}>0$, and $M$ linearly independent  solutions of \eqref{DiffEqn01} of the form
\begin{eqnarray}\label{tempsyst}
 e^{y_{j}(z,\epsilon;z_1)}, \quad j=1,\ldots, M,
\end{eqnarray}
 for  $ (z,\epsilon)\in  \Delta(z^{\prime},\eta^{\prime})\times \mathcal{G}(r,\alpha), \Delta(z^{\prime},\eta^{\prime})\subset  \Omega,j=1,\ldots ,M$ of \eqref{DiffEqn01} such that $\dsty y_j(z,\epsilon;z_1)=\frac{1}{\epsilon}\int_{z_1}^z\mathfrak{w}_j(t,\epsilon)dt$, with  $\mathfrak{w}_j(t,\epsilon)\in \mathcal{H}(\Delta(z^{\prime},\eta^{\prime})\times \mathcal{G}(r,\alpha))$ and $ \dsty
\mathfrak{w}_j(t,\epsilon)\sim \sum_{k=0}^{\infty}\mathfrak{b}_{j,k}(z)\epsilon^{k} $  as $\epsilon\rightarrow 0, \epsilon\in \mathcal{G}(r,\alpha)$, uniformly in  $z\in \Delta(z^{\prime},\eta^{\prime})$.   Here $z_1\in \Delta(z^{\prime},\eta^{\prime}) $ is an arbitrary point and  the path of integration in the above  integral is contained in  $\Delta(z^{\prime},\eta^{\prime})$.

The change  of variable  $v=e^{\frac{1}{\epsilon}\int^{z} (w(t,\epsilon)+ w_j(t))dt}$ transforms \eqref{DiffEqn01} into a generalized Ricatti equation of the form
\begin{multline}\label{RG}
\epsilon^{M-1}w^{(M-1)}=\\
 -\epsilon^{M-1} w_j^{(M-1)}-\epsilon^{M-1}\sum_{k=1}^{M-1}\frac{\rho_k }{\rho_M }(w+ w_j)^{(k-1)}- 
    \sum_{k=0}^{M-2}\epsilon^k\sum_{l=M-k}^{M}\frac{\rho_l }{\rho_M }{\bf B}_{l,M-k}(w+ w_j,\ldots,(w+ w_j)^{(l-M+k)})+ \frac{1}{\rho_M }.
\end{multline}

Next,   we transform \eqref{RG} into a non-linear first-order  system
\begin{eqnarray}\label{RGS}
\epsilon\frac{d\vec{g}}{dz}=\vec{G}_j(z,\vec{g},\epsilon). 
\end{eqnarray}
Further, we apply to \eqref{RGS}  a theorem, see below,  on the existence of asymptotic solutions \cite[Th.XII-1-2 p.374]{HSib99}  of Sibuya, see  also  \cite{Si58}.  Using  these solutions we obtain the WKB-expansion.

 \begin{theorem}[{\cite[Th.XII-1-2 p.374]{HSib99}, \cite{Si58}}]\label{SibuyaH}
 Consider the system
\begin{eqnarray}\label{SibSyst}
\epsilon v_j^{\prime}=f_j(z,v_1,\ldots,v_n,\epsilon), \quad j=1,\ldots,n,
\end{eqnarray}
where   $f_j(z,v_1,\ldots,v_n,\epsilon)$ are holomorphic functions of the  complex  variables $(z,v_1,\ldots,v_n,\epsilon)$ in the  domain
$$|z|<\delta_0, |\epsilon|<\rho_0, |\arg \epsilon|<\alpha_0,|v_j|<\gamma_0, \quad j=1,\ldots,n.$$
Set  
$$f_j(z,\vec{v},\epsilon)=f_{j0}(z,\epsilon) +\sum_{h=1}^{n}a_{jh}(z,\epsilon)v_h+\sum_{|p|\geq 2}f_{jp}(z,\epsilon)\vec{v}^p,$$
where $\vec{v}=(v_1,\ldots,v_n)\in \CC^n$, $p=(p_1,\ldots,p_n)$ is a multi--index of non--negative integers, and $\vec{v}^p=v_1^{p_1}\cdots v_n^{p_n}$ . Suppose that the following three conditions hold:
 \begin{itemize}
 \item[I)] Each function $f_j(z,\vec{v},\epsilon) (j=1,\ldots,n$) has an asymptotic expansion  
 $$f_j(z,\vec{v},\epsilon)\sim \sum_{\nu=0}^{\infty}\hat{f}_{j\nu}(z,\vec{v})\epsilon^{\nu},$$
in the sense of Poincar\'e as $\epsilon\rightarrow0$ in the sector 
\begin{eqnarray}\label{sector}
  |\epsilon|<\rho_0, |\arg \epsilon|<\alpha_0,
\end{eqnarray}
 where the coefficients $\hat{f}_{j\nu}(z,\vec{v})$ are holomorphic in the domain $|z|<\delta_0, |\vec{v}|<\gamma_0$.  Furthermore, we assume that  $\hat{f}_{j0}(z,\vec{0})=0,   j=1,\ldots, n, \quad  \mbox{for} \quad |z|<\delta_0$.

 \end{itemize}
 
\noindent  Let $A(z,\epsilon)$ be the $n\times n$ matrix whose $(j,k)$ entry is $a_{jk}(z,\epsilon)$. Then, $A(z,\epsilon)$ admits an asymptotic expansion,
 $$A(z,\epsilon)\sim \sum_{\nu=0}^{\infty}\epsilon^{\nu}A_{\nu}(z),$$
 as $\epsilon\rightarrow 0$ in the sector \eqref{sector}, where the entries of coefficients matrix $A_{\nu}(z)=(a_{jk\nu}(z))$ are holomorphic in the domain $|z|<\delta_0, |\vec{v}|<\gamma_0$.  
 
 \begin{itemize}
 
  \item[II)]  The eigenvalues of the  matrix $A_0(0)$  are such that  $\lambda_j\neq 0, j=1,\ldots, n$.

 \end{itemize}
 
Then,  system  \eqref{SibSyst} has a solution $v_j=p_j(z),  j=1,\ldots, n$ such that 
 \begin{itemize} 
 
  \item[i)]   $p_j(z,\epsilon)$ are holomorphic in a domain 
$$|z|<\delta, |\epsilon|<\rho, |\arg \epsilon|<\alpha,$$
where $\delta,\rho$, and $\alpha$ are suitable positive constants such that $0<\delta\leq \delta_0, 0<\rho\leq \rho_0$, and $0<\alpha\leq \alpha_0$.

  \item[ii)]  $p_j(z,\epsilon)$ admit asymptotic expansions 
  $$p_j(z,\epsilon)\sim \sum_{\nu=1}^{\infty}p_{j\nu}(z)\epsilon^{\nu}, \quad j=1,\ldots ,n$$
as $\epsilon\rightarrow 0$ in the sector $0<|\epsilon|<\rho, |\arg \epsilon|<\alpha$, where the coefficients $p_{\nu j}(z)$ are holomorphic in the domain $|z|<\delta$.
  \end{itemize}
 
 \end{theorem}
 
\begin{remark} Notice that by considering the transformation $z\mapsto z+z^{\prime}$, the preceding theorem also holds in any disk $|z-z^{\prime}|<\delta_0$, provided that  the hypotheses are true for $z$ varying in small neighborhood of the origin.
%\leavevmode 
%\item \cite[Th.XII-1-2 p.374]{HSib99} requires the  assumption  that the matrix $A_0(0)$ have a decomposition  
%$$A_0(0)=diag[\lambda_1,\ldots,\lambda_n]+\mathcal{N},$$
% where $\mathcal{N}$ is a lower triangular  nilpotent matrix. Such tdecompositon is always possible, by virtue of 
%\end{itemize}

\end{remark}

The application of Theorem \ref{SibuyaH}  to  system \eqref{RGS} requires an intermediate transformation.   In this case,  the matrix $A_0$ in the  hypothesis II of Sibuya's Theorem  is invertible in $\Omega$, this is  studied in Lemmas \ref{compa1}  and \ref{pre}.

In Section \ref{Sab} we apply the  previous results  to obtain the asymptotic behavior of $Q_n$ in compact subsets of $\Omega$ which  will be done in several steps. We will denote by $z_1, c_n$, and $n_0$ some  unimportant numerical constants  whose values may change from lemma  to lemma. Using  Proposition \ref{Li} we prove first in Lemmas \ref{SuperSt} and \ref{PreMain} that for every fixed point $z^{\prime}\in \Omega$ we can find $ \Delta(z^{\prime},\eta)\subset \Omega$ such that the exponential solution $e^{y_1}$ of the system \eqref{tempsyst} is dominant and there exist $ c_n\in\CC$, and  $\tilde{y}_1(z,\epsilon_n;z_1),  z_1\in  \Delta(z^{\prime},\eta)$ such that 
\begin{equation}\label{Qnrep}
Q_{n}(z)=c_{n}e^{\tilde{y}_1(z,\epsilon_n;z_1)}, \quad \forall n>n_{0},  
\end{equation}
where $\dsty \tilde{y}_1(z,\epsilon_n;z_1)=\frac{1}{\epsilon_n}\int_{z_1}^z\mathfrak{w}_1(t,\epsilon_n)dt$ and $\dsty 
\mathfrak{w}_1(t,\epsilon_n)\sim \sum_{k=0}^{\infty}\mathfrak{b}_{1,k}(z)\epsilon_n^{k} $  as $n\rightarrow \infty$ uniformly in $ z\in  \Delta(z^{\prime},\eta)$. Next, in Lemma \ref{PreMainext}  we 
prove that for any closed Jordan curve $\gamma$ enclosing $\supp \mu$,   there exist $z_1\in ext(\gamma), n_0$,   $y$, and a constant $c_n$  such that   
\begin{equation}\label{Qnrep1}
Q_{n}(z)=c_{n}e^{y(z,\epsilon_n;z_1)}, \quad \forall n>n_{0},\forall z\in  ext( \gamma),
\end{equation}
where  $y(z,\epsilon_n;z_1)=\frac{1}{\epsilon_n}\int_{z_1}^z\mathfrak{w}(t,\epsilon_n)dt, \mathfrak{w}(t,\epsilon_n)\in \mathcal{H}(ext(\gamma))$,   $\dsty
\mathfrak{w}(z,\epsilon_n)\sim \sum_{k=0}^{\infty}\mathfrak{b}_{k}(z)\epsilon_n^{k} $ as $n\rightarrow \infty$ when  $z$ varies in a  compact subset  $K\subset ext(\gamma)$.  Here the integration path $\Gamma$   is chosen so that $\Gamma\subset ext(\gamma)\setminus \tau $, $\dsty y(z_0,\epsilon_n;z_1)=\lim_{\begin{subarray}{c}
                           z\rightarrow z_0\\
                           z\in \tau^{+}
                          \end{subarray}}y(z,\epsilon_n;z_1)$,  when  $z_0\in \tau$, where $\tau $ is  a Jordan arc connecting $\infty$ and any point $\omega\in \supp\mu$ oriented so that $\omega$ is the endpoint. To prove this, we consider a covering  of $K$ of the form $\{\Delta(z^{\prime},\eta(z^{\prime}))\}_{z^{\prime}\in K}$,  where $\Delta(z^{\prime},\eta(z^{\prime}))$ is such that the representation \eqref{Qnrep} holds. Using   the Heine--Borel theorem we obtain a finite subcover of $K$ and  by analytic continuation of a function element $(\tilde{y}_1(z,\epsilon_n;z_1),\Delta(z^{\prime},\eta(z^{\prime})))$ we obtain the function $y$, which also can  be continued analytically to $ext(\gamma)$. 
                          
Using  the representation \eqref{Qnrep},  in Lemma \ref{previct} we prove  that  we can find a primitive of the function $\mathfrak{w}$ such that    for $n$ sufficiently large
\begin{equation}\label{Qnrep2}
Q_n(z)=e^{ \frac{1}{\epsilon_n}y(z,\epsilon_n)}, \forall z\in ext(\gamma), 
\end{equation}
where $\dsty
y(z,\epsilon_n)\sim \sum_{k=0}^{\infty}\Phi_{k}(z)\epsilon_n^{k}, $ as $n\rightarrow \infty$ when $z$ varies in compact subsets  $K\subset ext(\gamma)$. Here  $\Phi_k(z)$ is the primitive of the function $\mathfrak{b}_k$ in $ext(\gamma)\setminus \tau$   such that $\dsty \lim_{z\rightarrow\infty}\Phi_k(z)-h_k\ln z=0$, $\dsty \Phi_k(z)=\lim_{\begin{subarray}{c}
                           z\rightarrow z_0\\
                           z\in \tau^{+}
                          \end{subarray}}\Phi_k(z)$,  when  $z_0\in \tau\cap ext(\gamma)$, being  $\tau$  any Jordan arc  from  $-\infty$ to $\omega\in \supp \mu$ such that $]-\infty,p]\subset \tau$, for some $p\in \RR$ and    $h_k$ a constant  properly chosen.  The final step is the proof of Theorem \ref{Main} that follows by transforming \eqref{Qnrep2} appropriately.

In Section \ref{Stk} we discuss the Stokes phenomena   and in Section \ref{Ap} we present some applications.

\section{The existence of linearly independent exponential  solutions }\label{Li}

%Let $K\subset \Omega$ be  a compact set, for each $j=1,\dots,M$ pick a single value analytic branch of the function $\left(\sqrt[M]{\rho_M(z)}\right)^{-1}, z\in K$, define 
%\begin{eqnarray}
%w_j=\left(\sqrt[M]{\rho_M(z)}\right)^{-1}, z\in \Omega.
%\end{eqnarray}

\begin{lemma}\label{Fj}
For $M\geq 2$, let  $U\subset\Omega$ be an open subset, $S\subset \CC$, and   $w(t,\epsilon):U\times S\rightarrow \CC$  be a  holomorphic  function in the variable  $t$. Then   $v=e^{\frac{1}{\epsilon}\int^{z} (w(t,\epsilon)+ w_j(t))dt}$ is  a solution of 
\begin{eqnarray}\label{DF}
v^{(M)} +\sum_{k=1}^{M-1} \frac{\rho_k }{\rho_M } v^{(k)} - \frac{v }{\epsilon^{M} \rho_M  }=0,
\end{eqnarray}
if and only if, $w$ is a solution of 
$$\epsilon^{M-1}w^{(M-1)}=F_j(\epsilon,w,\ldots,w^{(M-2)}). $$
Here 
\begin{multline*}
F_j(\epsilon,w,\ldots,w^{(M-2)})=\\
 -\epsilon^{M-1} w_j^{(M-1)}-\epsilon^{M-1}\sum_{k=1}^{M-1}\frac{\rho_k }{\rho_M }(w+ w_j)^{(k-1)}- 
    \sum_{k=0}^{M-2}\epsilon^k\sum_{l=M-k}^{M}\frac{\rho_l }{\rho_M }{\bf B}_{l,M-k}(w+ w_j,\ldots,(w+ w_j)^{(l-M+k)})+ \frac{1}{\rho_M },
\end{multline*}
and $w_j$ is as in \eqref{wj}.

\end{lemma}
\proof
Let us denote by $P_{k}(y^{\prime},\dots ,y^{(k)})$ the polynomial in the   variables $(y^{\prime},\dots,y^{(k)})$  defined for all integers $1\leq k\leq M$ by the relation

$$ P_{k}(y^{\prime},\dots ,y^{(k)})={e^{-y}} \, \left(e^{y} \right)^{(k)}. $$

According to the Fa\`a di Bruno formula \cite[Th.A p.137]{comtet74}, $P_{k}(y^{\prime},\dots ,y^{(k)})$ can be expressed as
\begin{equation} \label{faa}
P_{k}(y^{\prime},\dots ,y^{(k)})=\sum_{l=1}^{k} {\bf B}_{k,l}(y^{\prime},\ldots,y^{(k-l+1)}).
\end{equation}

  The relation \eqref{faa}   and the change  of the variable  $v=e^{y}$ provide  that \eqref{DF} can be expressed as

\begin{equation} \label{DiffEqn03}
\sum_{k=1}^{M} \frac{\rho_k }{\rho_M } P_{k}(y^{\prime},\dots ,y^{(k)})- \frac{1}{\rho_M  \epsilon^{M}}=0.
\end{equation}

By multiplying the relation  \eqref{DiffEqn03} by $\epsilon^{M}$ and using the expression  for $P_k$, one has that  $v=e^y, 
y^{\prime} =(w+ w_j) \epsilon^{-1}$ is a solution to \eqref{DF} if and only if $w$ is a solution of the  generalized Riccati equation
\begin{eqnarray}\label{Operator01}
 \sum_{k=1}^{M} \epsilon^{M}\frac{\rho_k }{\rho_M } P_{k}(\epsilon^{-1}(w+ w_j),\dots ,\epsilon^{-1}(w+ w_j)^{(k-1)})- \frac{1}{\rho_M  }= 0,
\end{eqnarray}
where 
\begin{eqnarray*}
P_{k}(\epsilon^{-1}(w+ w_j),\dots ,\epsilon^{-1}(w+ w_j)^{(k-1)})=\sum_{l=1}^{k}\epsilon^{-l}{\bf B}_{k,l}(w+ w_j,\ldots,(w+ w_j)^{(k-l)}).
\end{eqnarray*}
By rearranging conveniently, we have that the relation \eqref{Operator01}  can be expressed equivalently as 
\begin{eqnarray*}
 \sum_{k=0}^{M-1}\epsilon^k\sum_{l=M-k}^{M}\frac{\rho_l }{\rho_M }{\bf B}_{l,M-k}(w+ w_j,\ldots,(w+ w_j)^{(l-M+k)})- \frac{1}{\rho_M  }=0,
 \end{eqnarray*}
that is, 
\begin{multline*}
\epsilon^{M-1}(w+ w_j)^{(M-1)}+\\
\epsilon^{M-1}\sum_{k=1}^{M-1}\frac{\rho_k }{\rho_M }(w+ w_j)^{(k-1)}+ 
    \sum_{k=0}^{M-2}\epsilon^k\sum_{l=M-k}^{M}\frac{\rho_l }{\rho_M }{\bf B}_{l,M-k}(w+ w_j,\ldots,(w+ w_j)^{(l-M+k)})- \frac{1}{\rho_M  }=0.
\end{multline*} 
\lqqd

\begin{lemma}\label{compa}
For $M\geq 3$, let $F_j; j=1,\ldots,M$ be defined as in Lemma \ref{Fj}.  Then, the $(M-1)$-st order differential equation 
\begin{eqnarray}\label{Feqj}
\epsilon^{M-1}w^{(M-1)}=F_j(\epsilon,w,\ldots,w^{(M-2)}), 
\end{eqnarray}
can be expressed in the companion form as 
\begin{eqnarray*}
\epsilon\frac{d\vec{g}}{dz}=\vec{G}_j(z,\vec{g},\epsilon),
\end{eqnarray*}
where $\vec{g}=(g_0,g_1,\ldots,g_{M-2})^t, \vec{G}_j(z,\vec{g},\epsilon)=(g_1+\epsilon w^{\prime}_j,g_2,\ldots,g_{M-2},h_j(z,\vec{g},\epsilon)),$ and 
\begin{multline}\label{hj}
h_j(z,\vec{g},\epsilon)= - \sum_{k=2}^{M} {\bf B}_{M,k}(g_0+w_j,g_1,\ldots,g_{M-k})+\frac{1}{\rho_M }-\sum_{k=1}^{M-2}\epsilon^{k}\left(\frac{\rho_{M-k}}{\rho_M} g_{M-k-1}+\right.\\
 \left. \sum_{l=k}^{M-2}\frac{\rho_{M-k} }{\rho_M }{\bf B}_{M-k,M-l}(g_0+w_j,g_1,\ldots,g_{l-k})   \right)-\epsilon^{M-1}\left(w_j^{(M-1)}+\frac{\rho_1}{\rho_M}(g_0+w_j)\right).   
\end{multline}
\end{lemma}
\proof
Consider the change of variables:
\begin{eqnarray*} 
\dsty g_0=w,g_1=\epsilon (g_0^{\prime}+w^{\prime}_j),g_2=\epsilon g_1^{\prime},\ldots
\end{eqnarray*}
then equation \eqref{Feqj} transform into 
$$\epsilon g^{\prime}_{M-2}=\dsty F_{j}\left(\epsilon, g_0+w_j,\frac{g_1}{\epsilon},\ldots,\frac{g_{M-2}}{\epsilon^{M-2}}\right).$$
Further,   using \cite[Th.A p.134]{comtet74} we have that 
 \begin{multline*}
F_{j}\left(\epsilon, g_0+w_j,\frac{g_1}{\epsilon},\ldots,\frac{g_{M-2}}{\epsilon^{M-2}}\right)=\\
 -\epsilon^{M-1} \left(w_j^{(M-1)}+\frac{\rho_1}{\rho_M}(g_0+w_j)\right)-\sum_{k=2}^{M-1}\frac{\rho_k }{\rho_M }\epsilon^{M-k}g_{k-1}-
 \sum_{k=0}^{M-2} \epsilon^k\left( \sum_{l=k}^{M-2}\frac{\rho_{M-k} }{\rho_M }{\bf B}_{M-k,M-l}(g_0+w_j,g_1,\ldots,g_{l-k})\right)+ \frac{1}{\rho_M }\\
=- \sum_{k=2}^{M} {\bf B}_{M,k}(g_0+w_j,g_1,\ldots,g_{M-k})+\frac{1}{\rho_M }-\sum_{k=1}^{M-2}\epsilon^{k}\left(\frac{\rho_{M-k}}{\rho_M} g_{M-k-1}+\right.\\
 \left. \sum_{l=k}^{M-2}\frac{\rho_{M-k} }{\rho_M }{\bf B}_{M-k,M-l}(g_0+w_j,g_1,\ldots,g_{l-k})   \right)-\epsilon^{M-1}\left(w_j^{(M-1)}+\frac{\rho_1}{\rho_M}(g_0+w_j)\right).
\end{multline*}
Thus,   the $(M-1)$-st order differential equation 
$$\epsilon^{M-1}w^{(M-1)}=F_j(\epsilon,w,\ldots,w^{(M-2)}), $$
can be expressed  as 
\begin{eqnarray*} 
\epsilon\frac{d\vec{g}}{dz}=
(g_1+\epsilon w^{\prime}_j,g_2,\ldots,g_{M-2},h_j(z,\vec{g},\epsilon))
\end{eqnarray*}
where $\vec{g}=(g_0,g_1,\ldots,g_{M-2})^t$ and $h_j$ is given by \eqref{hj}. \lqqd

The function $\vec{G}_j$ defined in Lemma \ref{compa} will play an important role;  for $M\geq 2$  we define  the function $\vec{G}_j$ as follows. Notice that  when $M=2$ the differential equation 
$$\epsilon^{M-1}w^{(M-1)}=F_j(\epsilon,w,\ldots,w^{(M-2)}), $$
by Lemma \ref{Fj} reduces to    the usual Riccati equation 
\begin{eqnarray}\label{m=2}
\dsty \epsilon w^{\prime}=-(w+w_j)^2+\frac{1}{\rho_2 }-\epsilon w_j^{\prime}-\frac{\rho_1 }{\rho_2 }\epsilon (w+w_j).
\end{eqnarray}

\begin{definition}\label{Gj}
For $M\geq 3$, define  $\vec{G}_j, 1\leq j\leq M$  as in  Lemma \ref{compa}.  For $M=2$, using \eqref{m=2},  we define  $\vec{G_j}=G_j$  as
\begin{eqnarray*}
\dsty G_j(z,g_0,\epsilon)=-(g_0+w_j)^2+\frac{1}{\rho_2 }-\epsilon w_j^{\prime}-\frac{\rho_1 }{\rho_2 }\epsilon (g_0+w_j), \quad j=1,2.
\end{eqnarray*}
\end{definition}

\begin{lemma}\label{compa1}
 Take  $\vec{G}_j,j=1,\ldots,M$ be as in Definition \ref{Gj}  and define $A_{j}(z,\epsilon)=\frac{\partial \vec{G}_j}{\partial \vec{g}}(z,\vec{0},\epsilon)$. Then,   $\dsty A_{j,0}(z)=\lim_{\epsilon\rightarrow 0}A_{j}(z,\epsilon)$ is   invertible in $ \Omega$.
 \end{lemma}
\proof
For $M=2$, it is immediate  from the Definition \ref{Gj} that $\dsty \lim_{\epsilon\rightarrow 0}\frac{\partial G_j}{\partial g_0}(z,0,\epsilon) \neq 0, z\in  \Omega,  j=1,2$.  Suppose now that $M\geq 3$, from   the expression for $\vec{G}_j$,  we get
\begin{eqnarray}\label{G0}
\dsty \frac{\partial \vec{G}_j}{\partial \vec{g}}(z,\vec{g},\epsilon)=\begin{pmatrix}
           0        &1           &     0       &   \cdots        &          0  \\
            0       & 0         &      1       &      \cdots      &          0   \\
   & &  \vdots  & &  \\
 \dsty   \frac{\partial h_j(z,\vec{g},\epsilon)}{\partial g_0}    & \dsty  \frac{\partial h_j(z,\vec{g},\epsilon)}{\partial g_1}   &  \cdots      &    & \dsty\frac{\partial h_j(z,\vec{g},\epsilon)}{\partial g_{M-2}}
\end{pmatrix}.
\end{eqnarray}

Using the definition of the complete Bell polynomials \cite[3b) p.134]{comtet74} we get the identity 
\begin{eqnarray}\label{Bellid}
\frac{\partial {\bf Y}_{n}}{\partial x_k}(x_1,\ldots, x_n)=\binom{n}{k} {\bf Y}_{n-k}(x_1,\ldots, x_{n-k}).
\end{eqnarray}
On the other hand, from the identity \cite[3c) p.134]{comtet74} we can write  $h_j$ as, 
$$h_j(z,\vec{g},\epsilon)=-{\bf Y}_{M}(g_0+w_j,g_1,\ldots, g_{M-2},0)+  \pi_{M-2}(\vec{g},\epsilon),$$
where   $\pi_{M-2}(\vec{g},\epsilon)$ is a polynomial such that $\dsty \frac{\partial \pi_{M-2}}{\partial g_k}(\vec{g},0)=0, 0\leq k\leq M-2$. Therefore,  by    \eqref{Bellid} 
$$  \frac{\partial h_j}{\partial g_k}(z,\vec{g},\epsilon)=-\binom{M}{k+1} {\bf Y}_{M-k-1}(g_0+w_j,\ldots, g_{M-k-1}) +\frac{\partial \pi_{M-2}(\vec{g},\epsilon)}{\partial g_k},  \quad 0\leq k\leq M-2.$$
By substituting this last expression  in  \eqref{G0} and taking the limit as $\epsilon \rightarrow 0$, we obtain 
\begin{eqnarray*}
\dsty A_{j,0}(z)=
\begin{pmatrix}
           0        &1           &     0       &   \cdots        &          0  \\
        0            & 0         &1            &      \cdots      &          0   \\
            &      &    \vdots           &        &     \\
     -\binom{M}{1}w_j^{M-1} &   -\binom{M}{2}w_j^{M-2}  &  \cdots      &  -\binom{M}{M-2}w_j^2& -\binom{M}{M-1}w_j
\end{pmatrix}.
\end{eqnarray*}

Now, to prove that $A_{j,0}$ is invertible in $ \Omega$ we prove that $0$ is not an eigenvalue.  Indeed,  the matrix $A_{j,0}$ is the companion matrix of the polynomial 
$$p(\lambda)=\lambda^{M-1}+\binom{M}{M-1}w_j\lambda^{M-2}+\ldots +\binom{M}{1}w_j^{M-1},$$
cf. \cite[Def 3.3.13 p.195]{hj}.  Hence, $p$ is the characteristic polynomial of the matrix $A_{j,0}$. On the other hand, we have that   $(\lambda+w_j)^{M}=\lambda p(\lambda)+w_j^M$.  A straightforward calculation shows that the roots of $p(\lambda)$ are given by the products $(\lambda_kw_j)_{k=1}^{M}$, where $\lambda_k=(\omega_k-1)$ and $\omega_k=e^{\frac{\imath 2k\pi}{M}}, k=1,\ldots,M-1$. Since $\lambda_kw_j(z)\neq 0, z\in \Omega$, the result is proven. \lqqd

\begin{lemma}\label{pre}
Let $\vec{G}_j,j=1,\ldots,M$ be as in Definition \ref{Gj} and  $z^{\prime}\in\Omega, \eta>0$ be   such that   $\Delta(z^{\prime},\eta)\subset   \Omega$. Then,   there exist   $\eta^{\prime},r,\alpha>0$ such that the system  of  differential equations
$$
\dsty \epsilon \frac{d\vec{g}}{dz}=\vec{G}_{j}(z,\vec{g},\epsilon), \quad z\in \Delta(z^{\prime},\eta)
$$
has a    solution    $\vec{g}_j$    holomorphic in $(z,\epsilon)\in  \Delta(z^{\prime},\eta^{\prime})\times \mathcal{G}(r,\alpha)$, for every  $1\leq  j\leq M$, admitting  a uniform  asymptotic expansion 
$$\vec{g}_j(z,\epsilon)\sim \sum_{k=1}^{\infty}\epsilon^k\vec{g}_{j,k}(z), \quad  \mbox{as  $\epsilon\rightarrow 0$ in $  \mathcal{G}(r,\alpha)$}, z\in \Delta(z^{\prime},\eta^{\prime}).$$

\end{lemma}

\proof 
 
  Consider first that $M\geq3$. By the Definition \ref{Gj},  $\vec{G}_j$ can be expressed as 
\begin{eqnarray}\label{writingG}
\vec{G}_j(z,\vec{g},\epsilon)=   \vec{f}_{j,0}(z,\epsilon) +A_{j,0}(z)\vec{g}+\vec{f_j}(z,\epsilon,g_0,\ldots,g_{M-2}),   \quad (z,\epsilon)\in \Delta(z^{\prime},\eta)\times \mathcal{G}(r^{\prime},\alpha^{\prime}),
\end{eqnarray} 
where $A_{j,0}(z)$ is defined as in   Lemma  \ref{compa1}, $\vec{f}_{j,0}(z,\epsilon)$ is a vector function whose components are polynomials in the variable $\epsilon$ with coefficients in $\mathcal{H}(\Delta(z^{\prime},\eta))$ such that  $\vec{f}_{j,0}(z,0)=\vec{0}$, $\vec{f_j}$  is a vector function whose components are  polynomials in the variables $\epsilon,g_0,\ldots,g_{M-2}$  with coefficients in $\mathcal{H}(\Delta(z^{\prime},\eta))$  such that $\vec{f_j}(z,0,g_0,\ldots,g_{M-2})=\vec{0}$.  We have  that for  each  $j=1,\ldots,M$ the right--hand side of \eqref{writingG}   satisfies the hypothesis  of   Theorem \ref{SibuyaH} in the open set  $\Delta(z^{\prime},\eta)\times \mathcal{G}(r^{\prime},\alpha^{\prime})$. Hence, for each fixed  $j=1,\ldots,M$, there exist $\eta_j\leq \eta, \alpha_j\leq \alpha^{\prime} , r_j\leq r^{\prime}$ such that   the  the differential equation \eqref{writingG}
has a    solution    $\vec{g}_j=(g_{0,j},\ldots, g_{{M-2},j})$    holomorphic in $(z,\epsilon)\in  \Delta(z^{\prime},\eta_j)\times \mathcal{G}(r_j,\alpha_j)$,  where   $g_{l,j}; 0 \leq l \leq M-2$ admits a uniform  asymptotic expansion 
$$g_{l,j}(z,\epsilon)\sim \sum_{k=1}^{\infty}\epsilon^kg_{l,j,k}(z), \quad  \mbox{as  $\epsilon\rightarrow 0$ in $  \mathcal{G}(r_j,\alpha_j)$}, z\in \Delta(z^{\prime},\eta_j),$$
hence, if  $\dsty r=\min_{j=1,\ldots,M}\{r_j\}, \dsty \alpha=\min_{j=1,\ldots,M}\{\alpha_j\}$ and $\dsty \eta^{\prime}=\min_{j=1,\ldots,M}\{\eta_j\}$ we obtain lemma for $M\geq 3$.

Finally, when $M=2$ one has by Definition \ref{Gj} that 
$$\epsilon \frac{dg_0}{dz}=-\epsilon w_j^{\prime}-\frac{\rho_1 }{\rho_2 }w_j\epsilon -2w_jg_0+\frac{\rho_1}{\rho_2 }g_0\epsilon-g_0^2, \quad   z\in (z,\epsilon)\in \Delta(z^{\prime},\eta)\times \mathcal{G}(r^{\prime},\alpha^{\prime}), j=1,2.$$
The conclusions of the lemma follows immediately again  from   Theorem \ref{SibuyaH}. \lqqd

\begin{proposition}\label{GL}
Take $M\geq 2$ and   $z^{\prime}\in \Omega$.  Then,  
\begin{itemize}

\item[a)]   there exist  $\alpha,r,\eta^{\prime}>0$, and $M$  solutions of \eqref{DiffEqn01} of the form   $ e^{y_{j}(z,\epsilon;z_1)}, (z,\epsilon)\in  \Delta(z^{\prime},\eta^{\prime})\times \mathcal{G}(r,\alpha), \Delta(z^{\prime},\eta^{\prime})\subset  \Omega,j=1,\ldots ,M$  such that $\dsty y_j(z,\epsilon;z_1)=\frac{1}{\epsilon}\int_{z_1}^z\mathfrak{w}_j(t,\epsilon)dt, \mathfrak{w}_j(t,\epsilon)\in \mathcal{H}(\Delta(z^{\prime},\eta^{\prime})\times \mathcal{G}(r,\alpha)),\,\,\, \dsty
\mathfrak{w}_j(t,\epsilon)\sim \sum_{k=0}^{\infty}\mathfrak{b}_{j,k}(t)\epsilon^{k} $  as $\epsilon\rightarrow 0, \epsilon\in \mathcal{G}(r,\alpha)$, uniformly when  $z$ varies in   $\Delta(z^{\prime},\eta^{\prime})$.   Here $z_1\in \Delta(z^{\prime},\eta^{\prime}) $ is an arbitrary point and  the integration  path in  the above integral is contained in  $\Delta(z^{\prime},\eta^{\prime})$;

\item[b)]
\begin{eqnarray*}%\label{bob2}  
\nonumber   \mathfrak{b}_{j,0}(z)&=& w_j(z), \quad  z\in \Delta(z^{\prime},\eta^{\prime}), \\
\nonumber \dsty \mathfrak{b}_{j,1}(z)&=&\frac{M-1}{2M}\frac{\rho^{\prime}_M(z)}{\rho_M(z)}-\frac{1}{M}\frac{\rho_{M-1}(z)}{\rho_{M}(z)}, \quad  z\in \Delta(z^{\prime},\eta^{\prime}),
\end{eqnarray*}
where $w_j$ is defined as in \eqref{wj};
\item[c)]  the $M$-tuple of functions $(e^{y_{j}(z,\epsilon;z_1)})_{j=1}^{M} $ is linearly independent  in $\Delta(z^{\prime},\eta^{\prime}) \times \mathcal{G}(r,\alpha)$ for $r$ small enough, where  $\{y_j\}_{j=1}^{M} $ is  as in item   a).

\end{itemize}
\end{proposition}
\proof

$a)$   By Lemma \ref{compa} and  relation \eqref{m=2}  we have that the $(M-1)$-st order differential equation
$$\epsilon^{M-1}w^{(M-1)}=F_j(\epsilon,w,\ldots,w^{(M-2)}), $$
 can be expressed in the companion form as 
$$\dsty \epsilon \frac{d\vec{g}}{dz}=\vec{G}_{j}(z,\vec{g},\epsilon),$$
where $\vec{G}_{j}$ is   as in Definition \ref{Gj}. Hence,    Lemma   \ref{pre} implies that   there exist    $\alpha,r,\eta^{\prime}>0$  and  a  holomorphic  function   $u_j$ in $\Delta(z^{\prime},\eta^{\prime})\times \mathcal{G}(r,\alpha),   \Delta(z^{\prime},\eta^{\prime})\subset \Omega$  such that 
\begin{eqnarray}\label{eqFsatisf}
\epsilon^{M-1}u_j^{(M-1)}=F_j(\epsilon,u_j,\ldots,u_j^{(M-2)}),  \quad 1\leq j\leq M,
\end{eqnarray}
where $u_j$ admits a uniform  asymptotic expansion 
\begin{eqnarray}\label{bmathfrak}
u_j(z,\epsilon)\sim \sum_{k=1}^{\infty}\epsilon^k\mathfrak{b}_{j,k}(z), \quad  \mbox{as  $\epsilon\rightarrow 0$ in $  \mathcal{G}(r,\alpha)$}, z\in \Delta(z^{\prime},\eta^{\prime}).
\end{eqnarray}
 Hence, by Lemma \ref{Fj} one has that if  
\begin{eqnarray}
\begin{aligned}\label{yjfraw} 
y_j(z,\epsilon;z_1)&=\epsilon^{-1}\int_{z_1}^{z} w_j(t) dt  +\epsilon^{-1}\int_{z_1}^{z}u_j (t,\epsilon)dt,\\
                    &\coloneqq \epsilon^{-1}\int_{z_1}^{z}\mathfrak{w}_j (t,\epsilon)dt\\
                    &=\sum_{k=0}^{\infty}\epsilon^{k-1}\int_{z_1}^{z}\mathfrak{b}_{j,k}(t)dt,
\end{aligned} 
\end{eqnarray}
  then  $ e^{y_{j}(z,\epsilon)}, (z,\epsilon)\in   \Delta(z^{\prime},\eta^{\prime})\times \mathcal{G}(r,\alpha)$ is a solution  of \eqref{DiffEqn01}, and $y_{j}$ admits a uniform  asymptotic expansion 
\begin{eqnarray}\label{yjfraw1}
y_{j}(z,\epsilon;z_1)\sim \dsty \sum_{k=0}^{\infty}\epsilon^{k-1}\int_{z_1}^{z}\mathfrak{b}_{j,k}(t)dt, \quad  \mbox{as  $\epsilon\rightarrow 0$ in $  \mathcal{G}(r,\alpha)$}, z\in \Delta(z^{\prime},\eta^{\prime}), 
\end{eqnarray} 
where  $\mathfrak{b}_{j,0}=w_j$.

$b)$
The first relation follows from \eqref{yjfraw} and \eqref{yjfraw1}. To prove the second assertion, notice that by   \eqref{eqFsatisf},  one has 
\begin{eqnarray}\label{derH0}
\frac{d}{d\epsilon}H_j(\epsilon,u_j,\ldots,u_j^{(M-1)})=0,
\end{eqnarray}
where $H_j(\epsilon,w,\ldots,w^{(M-1)})=\epsilon^{M-1}w^{(M-1)}-F_j(\epsilon,w,\ldots,w^{(M-2)})$. Further 
\begin{eqnarray}
\label{derH1} \frac{d}{d\epsilon}H_j(\epsilon,u_j,\ldots,u_j^{(M-1)})&= &\frac{\partial H_j}{\partial \epsilon}(\epsilon,u_j,\ldots,u_j^{(M-1)}) +\sum_{k=0}^{M-1} \frac{\partial H_j}{\partial w^{(k)}}(\epsilon,u_j,\ldots,u_j^{(M-1)}) \frac{d u_j^{(k)}}{d\epsilon} \\
\nonumber     &= & \frac{\rho_{M-1} }{\rho_M }{\bf B}_{M-1,M-1}(u_j+ w_j) +{\bf B}_{M,M-1}(u_j+ w_j,(u_j+ w_j)^{\prime})+
             I_j(\epsilon,u_j,\ldots,u_j^{(M-1)})\\
\nonumber      & & + \frac{\partial H_j}{\partial w}(\epsilon,u_j,\ldots,u_j^{(M-1)}) \frac{d u_j}{d\epsilon}   + \sum_{k=1}^{M-1} \frac{\partial H_j}{\partial w^{(k)}}(\epsilon,u_j,\ldots,u_j^{(M-1)}) \frac{du_j^{(k)}}{d\epsilon},               
\end{eqnarray}
where 
\begin{eqnarray*}
\nonumber I_j(\epsilon,u_j,\ldots,u_j^{(M-1)})&=& (M-1)\epsilon^{M-2} w_j^{(M-1)}+(M-1)\epsilon^{M-2}\sum_{k=1}^{M-1}\frac{\rho_k }{\rho_M }(u_j+ w_j)^{(k-1)}+ \\
\nonumber   & & \sum_{k=2}^{M-2}k\epsilon^{k-1}\sum_{l=M-k}^{M}\frac{\rho_l }{\rho_M }{\bf B}_{l,M-k}(u_j+ w_j,\ldots,(u_j+ w_j)^{(l-M+k)}),
\end{eqnarray*}
and $\dsty \frac{\partial H_j}{\partial w^{(k)}}$ denotes the partial derivative of $H_j$ in relation to the variable $w^{(k)}$, considering the  different derivatives of $w$ as independent functions.

On the other hand, using \cite[Th. XI--1--11 p. 347]{HSib99}, we  obtain the  relations,
\begin{eqnarray}
\begin{aligned}\label{derH2}
\lim_{\begin{subarray}{c}
                           \epsilon\rightarrow  0\\
                          \epsilon\in \mathcal{G}(r,\alpha)
                          \end{subarray}}   I_j(\epsilon,u_j,\ldots,u_j^{(M-1)})&= 0,  \\
\lim_{\begin{subarray}{c}
                           \epsilon\rightarrow  0\\
                          \epsilon\in \mathcal{G}(r,\alpha)
                          \end{subarray}}  \frac{\partial H_j}{\partial w}(\epsilon,u_j,\ldots,u_j^{(M-1)}) &= M(w_j)^{M-1},\\
\lim_{\begin{subarray}{c}
                           \epsilon\rightarrow  0\\
                          \epsilon\in \mathcal{G}(r,\alpha)
                          \end{subarray}} \sum_{k=1}^{M-1} \frac{\partial H_j}{\partial w^{(k)}}(\epsilon,u_j,\ldots,u_j^{(M-1)}) \frac{du_j^{(k)}}{d\epsilon}  &= 0.
\end{aligned}                          
\end{eqnarray}
 
By taking limits as  $\epsilon\rightarrow  0,  \epsilon\in \mathcal{G}(r,\alpha)$ in \eqref{derH1} and using  \eqref{bmathfrak},\eqref{derH0}, and \eqref{derH2}  one obtains 
$$0=  \frac{\rho_{M-1} }{\rho_M }{\bf B}_{M-1,M-1}(w_j) +{\bf B}_{M,M-1}( w_j,w_j^{\prime})+ M(w_j)^{M-1}\mathfrak{b}_{j,1}(z),$$ 
which  proves b).

$c)$ Using  b),  we see  that  the M-tuple  of functions $\{e^{y_{j}(z,\epsilon)}\}_{j=1}^{M} $  can be expressed as 
\begin{eqnarray*} 
 \left(e^{\frac{1}{\epsilon}\int_{z_1}^{z}w_1(t)dt+r_1(z,\epsilon)}, \ldots, e^{\frac{1}{\epsilon}\int_{z_1}^{z}w_M(t)dt+r_M(z,\epsilon)}\right),
 \end{eqnarray*} 
 where $r_j(z,\epsilon)=O(1)$ when $\epsilon\rightarrow 0, (z,\epsilon)\in \Delta(z^{\prime},\eta^{\prime})\times  \mathcal{G}(r,\alpha)$.   A straightforward calculation shows that the Wronskian of the system reads as 
\begin{multline*}
W\left(e^{\frac{1}{\epsilon}\int_{z_1}^{z}w_1(t)dt+r_1(z,\epsilon)}, \ldots, e^{\frac{1}{\epsilon}\int_{z_1}^{z}w_M(t)dt+r_M(z,\epsilon)}\right)=\prod_{k=1}^{M}e^{\frac{1}{\epsilon}\int_{z_1}^{z}w_k(t)dt+r_k(z,\epsilon)}\times \\
\left| \begin{array}{cccc}
           1                 &      \cdots          &                1  \\
           {\bf Y}_1\left(\frac{w_1(z)}{\epsilon}+r_1(z,\epsilon)\right)        &      \cdots      &          {\bf Y}_1\left(\frac{w_M(z)}{\epsilon}+r_M(z,\epsilon)\right)   \\
   & &  \vdots  &  \\
  {\bf Y}_{M-1}\left(\frac{w_1(z)}{\epsilon}+r_1(z,\epsilon),\ldots\right)        &  \cdots         &{\bf Y}_{M-1}\left(\frac{w_M(z)}{\epsilon}+r_M(z,\epsilon),\ldots\right)   
\end{array}\right|,
\end{multline*}
 using the identity \cite[(3c), p. 134]{comtet74}, one has  $\dsty  {\bf Y}_{n}(x_1,\ldots )= x_1^n+\sum_{k=1}^{n-1} {\bf B}_{n,k}(x_1,\ldots )$,  hence    the determinant $D$ in the above expression can be written as 
\begin{eqnarray*}
\begin{aligned}
D(z,\epsilon)&=\epsilon^{-\frac{M(M-1)}{2}}\left| \begin{array}{cccc}
           1        &1           &      \cdots          &                1  \\
       w_1(z)+o(\epsilon)        &  w_2(z)+o(\epsilon)          &      \cdots      &    w_M(z)+o(\epsilon)   \\
   & &  \vdots  &  \\
 \dsty   (w_1(z))^{M-1}+o(\epsilon)    & \dsty   (w_2(z))^{M-1}+o(\epsilon)       &  \cdots         &  (w_M(z))^{M-1}+o(\epsilon) 
\end{array}\right| \\
  &=\epsilon^{-\frac{M(M-1)}{2}}\left( V(z)   +o(\epsilon)\right),
\end{aligned}
\end{eqnarray*}
as $ \epsilon\rightarrow 0,   (z,\epsilon)\in \Delta(z^{\prime},\eta^{\prime})\times  \mathcal{G}(r,\alpha)$, where 
$$V(z)=\left| \begin{array}{cccc}
           1        &1           &      \cdots          &                1  \\
      w_1(z)     &  w_2(z)    &      \cdots      &       w_M(z)  \\
   & &  \vdots  &  \\
 \dsty   (w_1(z))^{M-1}   & \dsty   (w_2(z))^{M-1}     &  \cdots         &  (w_M(z))^{M-1} 
\end{array}\right|.$$
On the other hand,   for  the Vandermonde determinant $V$ we have that  $V(z)\neq 0, z\in \Delta(z^{\prime},\eta^{\prime})$. Therefore,  for $r$ small,  one has that $D\neq 0$ when  $(z,\epsilon)\in \Delta(z^{\prime},\eta^{\prime})\times  \mathcal{G}(r,\alpha)$.  Hence, 
$$W\left(e^{\frac{1}{\epsilon}\int_{z_1}^{z}w_1(t)dt+r_1(z,\epsilon)}, \ldots, e^{\frac{1}{\epsilon}\int_{z_1}^{z}w_M(t)dt+r_M(z,\epsilon)}\right) \neq 0, $$ 
 when  $(z,\epsilon)\in \Delta(z^{\prime},\eta^{\prime})\times  \mathcal{G}(r,\alpha)$ and  small $r$. 
This proves the assertion c).   \lqqd

%%%%%%%%%%%%%%%%%%%%%%%%%%%%%%%%%%%%%%%%%%%%%%%%%%%%%%%%%%%%%%%%%%%%%%%%%%%%%%%%%%%%%%%%%%%%%%%%%%%%%%%%%%%%%%%%%%%%%%%%%%%%%%%%%%%%%%%%%

%take a branch cut in $\CC$ defined by a Jordan arc connecting the zeros of the polynomial $\dsty \sum_{k=1}^{M}\rho_{k,k}(z)_k$ and   select 

\section {Asymptotic expansion}\label{Sab}

In this section we settle Theorem \ref{Main}, after proving  a series of preliminary lemmas. We   denote through this section the small parameter  
\begin{eqnarray}\label{epsilonn}
 \epsilon_n=\frac{1}{ \sqrt[M]{\lambda_{n}-\rho_{0,0}}},
\end{eqnarray}
where $\lambda_n$ is defined as in \eqref{lambn}; here  we  take the   branch of the root in \eqref{epsilonn} for which  $\epsilon_n$ coincides asymptotically with  $\dsty \frac{1}{n}$ when $n\rightarrow\infty$.  Recall   that $Q_n$ is a polynomial solution to \eqref{DiffEqn01} such that  $\lambda_{n}=\rho_{0,0}+\epsilon_n^{-M}$.

\begin{lemma}\label{SuperSt}
Let   $y_j(z,\epsilon;z_1), (z,\epsilon)\in \Delta(z^{\prime},\eta^{\prime})\times \mathcal{G}(r,\alpha), j=1,\ldots ,M$ be   as in item a) of  Proposition \ref{GL}. Then, 
\begin{itemize}
\item[a)]   there exist  open subsets $S_+(z_1),S_-(z_1)\subset \Delta(z^{\prime},\eta^{\prime})$ and $0<\alpha^{\prime}<\alpha$ such that 
 \begin{eqnarray}
\label{stokes1} \dsty \Re\left[\epsilon^{-1}\left(\int_{z_1}^{z}w_1(t)dt-\int_{z_1}^{z}w_j(t)dt\right)\right]>0, \quad \forall z\in  S_+(z_1),\epsilon\in \mathcal{G}(r,\alpha^{\prime}), j=2,\ldots, M,\\
\label{stokes2}  \dsty \Re\left[\epsilon^{-1}\left(\int_{z_1}^{z}w_1(t)dt-\int_{z_1}^{z}w_j(t)dt\right)\right]<0, \quad \forall z\in S_-(z_1), \epsilon\in \mathcal{G}(r,\alpha^{\prime}), j=2,\ldots, M;
\end{eqnarray}

\item[b)]  there exists  $z_1\in \Delta(z^{\prime},\eta^{\prime})$ such that $z^{\prime}\in S_+(z_1)$, where $S_+(z_1)$ is defined as in item a); 

\item[c)]  there exists $z^{\prime\prime}\in S_-(z_1)\setminus \supp\mu$  such that $Q^{(j)}_n(z^{\prime\prime})\neq 0$ for all $n\geq0$ and $0\leq j\leq M-2$; 

\item[d)]  
$
\dsty 
W(\dots,e^{y_{j-1}(z,\epsilon_n;z_1)},Q_n(z),e^{y_{j+1}(z,\epsilon_n;z_1)},\dots)= 
Q_n(z)\epsilon_n^{-\frac{M(M-1)}{2}}\left(H_j(z,\epsilon_n)   +o(\epsilon_n)\right)\prod_{\begin{subarray}{c}
                           k=1\\
                           k\neq j
                          \end{subarray}}^{M}e^{\frac{1}{\epsilon_n}\int_{z_1}^{z}w_k(t)dt+o(\epsilon_n)}$,\\
 where  
\begin{eqnarray}\label{H_j}
\dsty  H_j(z,\epsilon_n)=\left| \begin{array}{cccccc}
         \cdots  &    1        &1         &     1         &     \cdots   \\
     \cdots  &   w_{j-1}(z) &\dsty  \frac{\epsilon_n Q^{\prime}_n(z)}{Q_n(z)}   &    w_{j+1}(z) &     \cdots     \\
   & &  \vdots  &  \\
 \cdots  &   \dsty   (w_{j-1}(z))^{M-1}  & \dsty \frac{\epsilon_n^{M-1}Q^{(M-1)}_n(z)}{ Q_n(z)}    &  (w_{j+1}(z))^{M-1}  &     \cdots  
\end{array}\right|;
\end{eqnarray}

\item[e)]   
$\dsty  W(Q_n(z^{\prime\prime}),e^{y_{2}(z^{\prime\prime},\epsilon_n;z_1)},e^{y_{3}(z^{\prime\prime},\epsilon_n;z_1)},\dots) \neq 0$,  for all  $n$  large enough, 
where $z^{\prime\prime}$ is  as in item c);

\item[f)]  $
\dsty \lim_{n\rightarrow\infty}\left|\frac{W(\dots,e^{y_{j-1}(z^{\prime\prime},\epsilon_n;z_1)},Q_n(z^{\prime\prime}),e^{y_{j+1}(z^{\prime\prime},\epsilon_n;z_1)},\dots)}{W(Q_n(z^{\prime\prime}),e^{y_{2}(z^{\prime\prime},\epsilon_n;z_1)},e^{y_{3}(z^{\prime\prime},\epsilon_n;z_1)},\dots)}\right|=0, \quad j>1.
$

\end{itemize}
\end{lemma} 

\proof
a)  Since  $b_1(z)=\int_{z_1}^{z}w_1(t)dt$ is conformal in $\Delta(z^{\prime},\eta^{\prime})$ we have  that $U \equiv b_{1}(\Delta(z^{\prime},\eta^{\prime}))$ is an open set containing the origin.  Hence,  if $t=b_{1}(z)$
\begin{eqnarray}\label{preineq0}
\Re\left[\left(1-\exp\left(\imath \frac{2(j-1)\pi}{M}\right)\right)t\right]= \Re\left[\int_{z_1}^{z}w_1(t)dt-\int_{z_1}^{z}w_j(t)dt\right], \quad \forall z\in \Delta(z^{\prime},\eta^{\prime}), j=2,\ldots,M.
\end{eqnarray}

Notice  that   the sets
$$\{t\in U :\Re\left[\left(1-\exp\left(\imath \frac{2(j-1)\pi}{M}\right)\right)t\right]>0, j=2,\ldots,M\}$$
and 
 $$ \{t\in U :\Re\left[\left(1-\exp\left(\imath \frac{2(j-1)\pi}{M}\right)\right)t\right]<0, j=2,\ldots,M\}$$
are nonempty.  Therefore, for  $0<\alpha^{\prime}\leq \alpha$  small, one has that if $\epsilon\in \mathcal{G}(r,\alpha^{\prime})$, then 
\begin{eqnarray*}
\begin{aligned}
\mathcal{T}_+(z_1)=\{t\in U :\Re\left[\epsilon^{-1}\left(1-\exp\left(\imath \frac{2(j-1)\pi}{M}\right)\right)t\right]>0, j=2,\ldots,M\},\\
\mathcal{T}_-(z_1)=\{t\in U : \Re\left[\epsilon^{-1}\left(1-\exp\left(\imath \frac{2(j-1)\pi}{M}\right)\right)t\right]<0, j=2,\ldots,M\}
\end{aligned}
\end{eqnarray*}
are also  nonempty sets. Writing  $S_+(z_1)=b^{-1}_{1}(\mathcal{T}_+(z_1)), S_-(z_1)=b^{-1}_{1}(\mathcal{T}_-(z_1))$  and using  \eqref{preineq0} one obtains a).

b)   By \eqref{stokes1}  and \eqref{stokes2}, we see that $z_1$ is in the common boundary of the open sets $S_+(z_1)$ and $S_-(z_1)$. Hence, by taking   $z_1$   sufficiently close to $z^{\prime}$  we have   $z^{\prime}\in S_+(z_1)$.

c)  Let $\mathcal{A}$  be the set of zeros of $Q^{(j)}_n$, for all  $n\geq 0$ and $0\leq j\leq M-2$. Since $\mathcal{A}$  is at most countable we have that  $\mathcal{A}$ is of measure $0$. Since $\supp \mu$ is also a set of two-dimensional Lebesgue measure $0$, we can find $z^{\prime\prime}\in  S_-(z_1)$ such that $Q^{(j)}_n(z^{\prime\prime})\neq 0$, for all  $n\geq 0$ and $0\leq j\leq M-2$.

d)  Notice that  for $n$ large, $\epsilon_n\in  \mathcal{G}(r,\alpha)$. We have that 
\begin{multline*}
W(\dots,e^{y_{j-1}(z,\epsilon_n;z_1)},Q_n(z),e^{y_{j+1}(z,\epsilon_n;z_1)},\dots)= Q_n(z) \prod_{\begin{subarray}{c}
                           k=1\\
                           k\neq j
                          \end{subarray}}^{M}e^{\frac{1}{\epsilon_n}\int_{z_1}^{z}w_k(t)dt+r_k(z,\epsilon_n)}\times \\
\left| \begin{array}{cccccc}
     \cdots  &    1        & 1       &     1         &     \cdots  \\
       \cdots  &  {\bf Y}_1\left(\frac{w_{j-1}(z)}{\epsilon_n}+r_{j-1}(z,\epsilon_n)\right)       & \dsty  \frac{Q^{\prime}_n(z)}{ Q_n(z) }  &              {\bf Y}_1\left(\frac{w_{j+1}(z)}{\epsilon_n}+r_{j+1}(z,\epsilon_n)\right) & \cdots   \\
   & &  \vdots  &  \\
 \dsty   \cdots  &   {\bf Y}_{M-1}\left(\frac{w_{j-1}(z)}{\epsilon_n}+r_{j-1}(z,\epsilon_n),\ldots\right)    & \dsty \frac{Q^{(M-1)}_n(z)}{ Q_n(z) }    &{\bf Y}_{M-1}\left(\frac{w_{j+1}(z)}{\epsilon_n}+r_{j+1}(z,\epsilon_n),\ldots\right)   & \cdots  
\end{array}\right|.
\end{multline*}
Using the identity \cite[(3c), p. 134]{comtet74}, ones has  $\dsty  {\bf Y}_{n}(x_1,\ldots )= x_1^n+\sum_{k=1}^{n-1} {\bf B}_{n,k}(x_1,\ldots )$, hence
\begin{multline*}
W(\dots,e^{y_{j-1}(z,\epsilon_n;z_1)},Q_n(z),e^{y_{j+1}(z,\epsilon_n;z_1)},\dots)=
\epsilon_n^{-\frac{M(M-1)}{2}}Q_n(z) \prod_{\begin{subarray}{c}
                           k=1\\
                           k\neq j
                          \end{subarray}}^{M}e^{\frac{1}{\epsilon_n}\int_{z_1}^{z}w_k(t)dt+o(\epsilon_n)}\times \\
\left| \begin{array}{cccccc}
         \cdots  &    1        &   1    &     1         &     \cdots   \\
     \cdots  &    w_{j-1}(z)+o(\epsilon_n)  & \dsty \frac{\epsilon_nQ^{\prime}_n(z)}{Q_n(z)  }   &    w_{j+1}(z)+o(\epsilon_n) &     \cdots     \\
   & &  \vdots  &  \\
 \cdots  &   \dsty   (w_{j-1}(z))^{M-1}+o(\epsilon_n)    &\dsty \frac{\epsilon_n^{M-1} Q^{(M-1)}_n(z)}{ Q_n(z)  }    &  (w_{j+1}(z))^{M-1}+o(\epsilon_n)  &     \cdots  
\end{array}\right| \\
  =\epsilon_n^{-\frac{M(M-1)}{2}}Q_n(z)\prod_{\begin{subarray}{c}
                           k=1\\
                           k\neq j
                          \end{subarray}}^{M}e^{\frac{1}{\epsilon_n}\int_{z_1}^{z}w_k(t)dt+o(\epsilon_n)}\left( H_j(z,\epsilon_n)   +o(\epsilon_n)\right).
\end{multline*}
where $H_j$ is given by \eqref{H_j}.

e) Let $z^{\prime\prime}$ be as in c) and $H_j(z,\epsilon_n)$ as in  item d).  By \cite[Lem.10 \& Cor.4]{BerRull02}, the root counting measure of $Q_n^{(j)}, 0\leq j\leq M$ converges weakly to a measure  $\mu$ whose Cauchy transform equals to $\dsty \frac{1}{\sqrt[M]{\rho_{M}(z)}}$. Hence, for every $ 1\leq j\leq M-1 $ fixed, 
\begin{equation}
\begin{aligned}\label{tanj}
\dsty  \lim_{n\rightarrow \infty}\frac{ \epsilon_n^j Q_{n}^{(j)}(z^{\prime\prime})}{Q_{n}(z^{\prime\prime})}&=\left(\frac{1}{\sqrt[M]{\rho_{M}(z^{\prime\prime})}}\right)^j\\
                &= (w_1(z^{\prime\prime}))^j,
\end{aligned}
\end{equation} 
therefore,
$$\lim_{n\rightarrow\infty}H_1(z^{\prime\prime},\epsilon_n)= H_1(z^{\prime\prime}),$$
where 
\begin{eqnarray}\label{H1}
\dsty H_1(z^{\prime\prime})=\left| \begin{array}{cccccc}
1      &     \cdots  &    1        &     1         &     \cdots   \\
\dsty  w_1(z^{\prime\prime}) &   \cdots  &    w_j(z^{\prime\prime}) &  w_{j+1}(z^{\prime\prime}) &     \cdots     \\
   & &  \vdots  &  \\
\dsty  (w_1(z^{\prime\prime}))^{M-1}   &  \cdots  &   \dsty   (w_j(z^{\prime\prime}))^{M-1}  &  (w_{j+1}(z^{\prime\prime}))^{M-1}  &     \cdots  
\end{array}\right|.
\end{eqnarray}

Since   $H_1(z^{\prime\prime}) \neq 0$, item d)  implies  that  $\dsty W(Q_n(z^{\prime\prime}),e^{y_{2}(z^{\prime\prime},\epsilon_n;z_1)},e^{y_{3}(z^{\prime\prime},\epsilon_n;z_1)},\dots) \neq 0$ for all $n$ large enough.

f)
By item d) 
\begin{eqnarray}\label{prelim}
\left|\frac{W(\dots,e^{y_{j-1}(z^{\prime\prime},\epsilon_n)},Q_n(z^{\prime\prime}),e^{y_{j+1}(z^{\prime\prime},\epsilon_n)},\dots)}{W(Q_n(z^{\prime\prime}),\dots,e^{y_{j}(z^{\prime\prime},\epsilon_n)},e^{y_{j+1}(z^{\prime\prime},\epsilon_n)},\dots)}\right|=\left| \frac{H_j(z^{\prime\prime},\epsilon_n)+o(\epsilon_n)}{H_1(z^{\prime\prime},\epsilon_n)+o(\epsilon_n)}    \right|e^{\Re[\epsilon^{-1}_n(\int_{z_1}^{z^{\prime\prime}}w_1(t)dt-\int_{z_1}^{z^{\prime\prime}}w_j(t)dt)]}|1+o(\epsilon_n)|.
\end{eqnarray}

By  \eqref{H_j} and \eqref{tanj} we get 
\begin{eqnarray}\label{lim}
\lim_{n\rightarrow\infty}H_j(z^{\prime\prime},\epsilon_n)=0, \quad j>1.
\end{eqnarray}

On the other hand, for $n$ large enough, $\epsilon_n\in  \mathcal{G}(r,\alpha^{\prime})$. Since $z^{\prime\prime}\in S_-(z_1)$, formula  \eqref{stokes2} implies
\begin{eqnarray}\label{lim2}
\lim_{n\rightarrow\infty}e^{\Re[\epsilon^{-1}_n(\int_{z_1}^{z^{\prime\prime}}w_1(t)dt-\int_{z_1}^{z^{\prime\prime}}w_j(t)dt)]}=0.
\end{eqnarray}

Finally, by   \eqref{H1},\eqref{prelim},\eqref{lim}, and \eqref{lim2} we obtain f). \lqqd

\begin{lemma}\label{PreMain}
Take   $z^{\prime}\in \Omega$.  Then,   there exist  $ n_{0}\in \NN, 0<\eta, z_1\in \Delta(z^{\prime},\eta)$, and   a unique constant $c_{n}$ such that
\begin{equation*}
Q_{n}(z)=c_{n}e^{\tilde{y}_1(z,\epsilon_n;z_1)}, \quad \forall n>n_{0}, \quad  \forall z\in  \Delta(z^{\prime},\eta), 
\end{equation*}
  where $\tilde{y}_1(z,\epsilon_n;z_1)=\frac{1}{\epsilon_n}\int_{z_1}^z\tilde{\mathfrak{w}}_1(t,\epsilon_n)dt, \tilde{\mathfrak{w}}_1(t,\epsilon_n)\in \mathcal{H}(\Delta(z^{\prime},\eta))$, $\dsty
\tilde{\mathfrak{w}}_1(t,\epsilon_n)\sim \sum_{k=0}^{\infty}\mathfrak{b}_{1,k}(z)\epsilon_n^{k} $  as $n \rightarrow \infty$, uniformly in  $z\in  \Delta(z^{\prime},\eta)$.
\end{lemma}

\proof Let $\Delta(z^{\prime},\eta^{\prime})$ be as in item  a) of  Proposition \ref{GL}.  We have  that for $n$ large, $\epsilon_n\in  \mathcal{G}(r,\alpha)$, hence by items a) and c) of  Proposition \ref{GL},   there exist    $M$ linearly independent solutions of \eqref{DiffEqn01}  in   $\Delta(z^{\prime},\eta^{\prime})$,  $\{e^{y_1(z,\epsilon_n;z_1)},\ldots, e^{y_M(z,\epsilon_n;z_1)}\}$; and  $M$ unique constants   $c_{1}(n),\dots,c_{M}(n)$ such that
\begin{eqnarray}\label{lincomb}
Q_{n}(z)=c_{1}(n)e^{y_1(z,\epsilon_n;z_1)}+\dots +c_{M}(n)e^{y_M(z,\epsilon_n;z_1)},\quad    z\in  \Delta(z^{\prime},\eta^{\prime}),  n>n_0,
\end{eqnarray}
 where  
$y_j(z,\epsilon_n;z_1)=\frac{1}{\epsilon_n}\int_{z_1}^z \mathfrak{w}_j(t,\epsilon_n)dt,  \mathfrak{w}_j(t,\epsilon_n)\in \mathcal{H}(\Delta(z^{\prime},\eta^{\prime} )\times \mathcal{G}(r,\alpha))$, $\dsty
 \mathfrak{w}_j(t,\epsilon_n)\sim \sum_{k=0}^{\infty}\mathfrak{b}_{j,k}(z)\epsilon_n^{k} $  as $n \rightarrow \infty$, and $z_1 \in \Delta(z^{\prime},\eta^{\prime})$.

\begin{figure}[H]
\centering
        \includegraphics[width=0.4\textwidth]{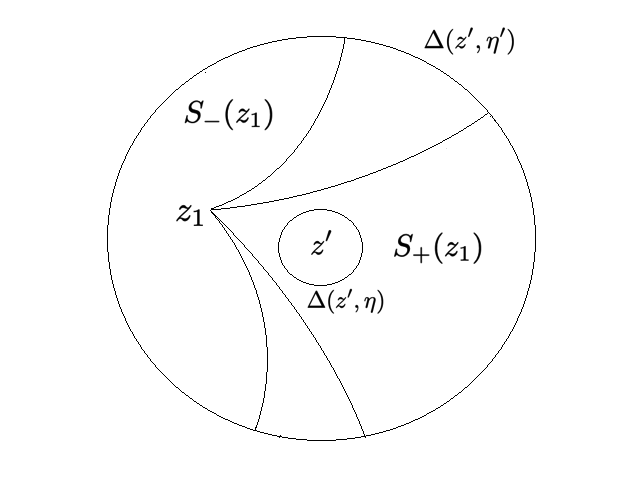}
         \caption{ Regions $\Delta(z^{\prime},\eta),  S_+(z_1), S_-(z_1)$}\label{fi00}
\end{figure}

  By item  b) of Lemma \ref{SuperSt} we can take   $z_1\in \Delta(z^{\prime},\eta^{\prime})$ such that $z^{\prime}\in S_+(z_1)$, where $S_+(z_1)$ is defined as in item a) of Lemma \ref{SuperSt}.  Choose $0<\eta<\eta^{\prime}$ so that $\Delta(z^{\prime},\eta)\subset S_+(z_1)$, and let $z^{\prime\prime}$ be  as in item  e) of Lemma \ref{SuperSt},  see Figure \ref{fi00}. Using Cramer's rule we have, 
\begin{equation*}
c_{j}(n)=\frac{W(\dots,e^{y_{j-1}(z^{\prime\prime},\epsilon_n;z_1)},Q_n(z),e^{y_{j+1}(z^{\prime\prime},\epsilon_n;z_1)},\dots)}{W( e^{y_{1}(z^{\prime\prime},\epsilon_n;z_1)}, \dots ,e^{y_{M}(z^{\prime\prime},\epsilon_n;z_1)})},
\end{equation*}
then, for $n$ large we can write  \eqref{lincomb}   as
\begin{multline}\label{Qna}
\dsty Q_{n}(z)=\\
c_{1}({n})e^{y_{1}(z,\epsilon_{n};z_1)}
 \left(\dsty 1+  \sum_{j=2}^{M}\frac{W(\dots,e^{y_{j-1}(z^{\prime\prime},\epsilon_n;z_1)},Q_n(z^{\prime\prime}),e^{y_{j+1}(z^{\prime\prime},\epsilon_n;z_1)},\dots)}{W(Q_n(z^{\prime\prime}),e^{y_{2}(z^{\prime\prime},\epsilon_n;z_1)},e^{y_{3}(z^{\prime\prime},\epsilon_n;z_1)},\dots)}e^{ y_j(z,\epsilon_{n};z_1)-y_1(z,\epsilon_{n};z_1)} \right), \quad z\in  \Delta(z^{\prime},\eta).
\end{multline}
On the other hand, since $z\in  \Delta(z^{\prime},\eta)\subset  S_+(z_1)$, we obtain from  the relation \eqref{stokes1} and item  f) of Lemma \ref{SuperSt} that the term inside the parenthesis in   \eqref{Qna} converges uniformly to $1$. Hence, for some $n_0$,
\begin{equation*} 
Q_{n}(z)=c^{\prime}_{n}e^{\tilde{y}_1(z,\epsilon_n;z_1)}, \quad \forall n>n_{0}, \quad  \forall z\in  \Delta(z^{\prime},\eta), 
\end{equation*}
 where $\tilde{y}_1(z,\epsilon_n;z_1)=\frac{1}{\epsilon_n}\int_{z_1}^z\tilde{\mathfrak{w}}_1(t,\epsilon_n)dt, \tilde{\mathfrak{w}}_1(t,\epsilon_n)\in \mathcal{H}(\Delta(z^{\prime},\eta))$, $\dsty
\tilde{\mathfrak{w}}_1(t,\epsilon_n)\sim \sum_{k=0}^{\infty}\mathfrak{b}_{1,k}(z)\epsilon_n^{k} $  as $n \rightarrow \infty$, and $z_1\in \Delta(z^{\prime},\eta^{\prime})$. Finally, we take the constant $z_1$ in the expression for $\tilde{y}_1$ so that $z_1\in \Delta(z^{\prime},\eta)$, which completes the proof of the lemma. \lqqd

 \begin{lemma}\label{compact}
Let $K\subset \Omega$ be a compact set.  Then, there exists a finite covering    $\{\Delta(z^{\prime}_i,\eta_i)\}_{i\in I}$ by  open disks such that  $\dsty K \subset \bigcup_{i\in I}\Delta(z^{\prime}_i,\eta_i)$  and $n_0$ such that 
$$
Q_{n}(z)=c_{n}(i)e^{\tilde{y}_{1,i}(z,\epsilon_n;z_{1,i})}, \quad \forall n>n_{0}, \forall z\in \Delta(z^{\prime}_i,\eta_i), i\in I,
$$
where $\tilde{y}_{1,i}(z,\epsilon_n;z_{1,i})=\frac{1}{\epsilon_n}\int_{z_{1,i}}^z\tilde{\mathfrak{w}}_{1,i}(t,\epsilon_n)dt$, $\tilde{\mathfrak{w}}_{1,i}(t,\epsilon_n)\in \mathcal{H}(\Delta(z^{\prime}_i,\eta_i))$,    $\dsty
\tilde{\mathfrak{w}}_{1,i}(t,\epsilon_n)\sim \sum_{k=0}^{\infty}\mathfrak{b}_{1,i,k}(z)\epsilon_n^{k} $  as $n\rightarrow \infty$ for $z$ varying  in compact subsets of $\Delta(z^{\prime}_i,\eta_i)$, and $z_{1,i}\in \Delta(z^{\prime}_i,\eta_i)$.
\end{lemma}

\proof
By Lemma \ref{PreMain}, for each $z^{\prime}\in \Omega $, there exists $\eta>0, \tilde{y}_1$, and $n_0(z^{\prime})$ such that 
\begin{equation*}
Q_{n}(z)=c_{n}e^{\tilde{y}_1(z,\epsilon_n;z^{\prime}_{1})}, \quad \forall n>n_{0}(z^{\prime}), \quad  \forall z\in  \Delta(z^{\prime},\eta(z^{\prime})), 
\end{equation*}
  where $\tilde{y}_1(z,\epsilon_n;z_1)=\frac{1}{\epsilon_n}\int_{z_1}^z\tilde{\mathfrak{w}}_1(t,\epsilon_n)dt$, $\tilde{\mathfrak{w}}_1(t,\epsilon_n)\in \mathcal{H}(\Delta(z^{\prime},\eta(z^{\prime})))$, $\dsty
\tilde{\mathfrak{w}}_1(t,\epsilon_n)\sim \sum_{k=0}^{\infty}\mathfrak{b}_{1,k}(z)\epsilon_n^{k} $  as $n \rightarrow \infty$ for $z$ varying  in compact subsets of $\Delta(z^{\prime},\eta(z^{\prime}))$, and $z_1\in \Delta(z^{\prime},\eta(z^{\prime}))$.

Let $U$ be an open set such that $K\subset U$ and consider the   open  covering  $\dsty K\subset \bigcup_{z^{\prime}\in U}\Delta(z^{\prime},\eta(z^{\prime}))$.  By   the Heine-Borel theorem we can find  a finite  open subcover such that $K\subset \bigcup\limits_{i\in I}\Delta(z_i^{\prime},\eta_i)$, where   $I$  is a finite set of indices.  By taking $n_0=\max\limits_{i\in I}n_0(z^{\prime}_i)$ we complete the proof of the lemma.  \lqqd

\begin{lemma}\label{PreMainext}
Let     $\gamma$ be a closed Jordan curve enclosing $\supp \mu$  and $\Omega_{\tau}=\Omega\setminus \tau$;  where $\tau $  is any Jordan arc connecting $\infty$ and any point $\omega\in \supp\mu$  and oriented so that $\omega$ is its endpoint. Then, 
\begin{itemize}
\item[a)]  there exist $z_1\in ext(\gamma), n_0, y$, and a constant $c_n$  such that   
\begin{equation*}
Q_{n}(z)=c_{n}e^{y(z,\epsilon_n;z_1)}, \quad \forall n>n_{0},\forall z\in  ext( \gamma),
\end{equation*}
where  $y(z,\epsilon_n;z_1)=\frac{1}{\epsilon_n}\int_{z_1}^z\mathfrak{w}(t,\epsilon_n)dt$, $\mathfrak{w}(t,\epsilon_n)\in \mathcal{H}(ext(\gamma))$,   $\dsty
\mathfrak{w}(z,\epsilon_n)\sim \sum_{k=0}^{\infty}\mathfrak{b}_{k}(z)\epsilon_n^{k} $ as $n\rightarrow \infty$ when $z$ varies in compact subsets  $K\subset ext(\gamma)$.  Here the integration  path   of the above integral is contained in  $ ext(\gamma)\setminus \tau $, and $\dsty y(z_0,\epsilon_n;z_1)=\lim_{\begin{subarray}{c}
                           z\rightarrow z_0\\
                           z\in \tau^{+}
                          \end{subarray}}y(z,\epsilon_n;z_1)$,  when  $z_0\in \tau$.  
                          
 \item[b)] further,
 \begin{eqnarray*}%\label{bob2}  
 \nonumber \mathfrak{b}_0(z)&=&w_1(z),  \quad z\in ext(\gamma),\\
\nonumber  \mathfrak{b}_1(z)&=& \left(\frac{M-1}{2M}\frac{\rho^{\prime}_M(z)}{\rho_M(z)}-\frac{1}{M}\frac{\rho_{M-1}(z)}{\rho_{M}(z)}\right), \quad z\in ext(\gamma).
\end{eqnarray*} 

\end{itemize}

\end{lemma}
\proof
a) By \cite[Lem. 9]{BerRull02}, we can find a Jordan closed curve $\gamma_1$ enclosing $\gamma$ and such   that for  $n>n_0$ the zeros of $Q_n$ are contained in $int(\gamma_1)$, define the compact set $K_0=\overline{int(\gamma_1)\setminus int(\gamma)}$. By Lemma \ref{compact} there exists a finite covering    $\mathfrak{S}=\{\Delta(z^{\prime}_i,\eta_i)\}_{i\in I_0}$ by open disks 
 such that   $\dsty K_0 \subset \bigcup_{i\in I_0}\Delta(z^{\prime}_i,\eta_i)$  and 
$$
Q_{n}(z)=c_{n}(i)e^{y_{1,i}(z,\epsilon_n;z_{1,i})}, \quad \forall n>n_{0}(i), \forall z\in  \Delta(z^{\prime}_i,\eta_i), i\in I_0,
$$
where 
\begin{eqnarray}\label{equa30}  
y_{1,i}(z,\epsilon_n;z_{1,i})=\frac{1}{\epsilon_n}\int_{z_{1,i}}^z\mathfrak{w}_{1,i}(t,\epsilon_n)dt, \quad   z,z_{1,i}\in \Delta(z^{\prime}_i,\eta_i), \mathfrak{w}_{1,i}(t,\epsilon_n)\in \mathcal{H}(\Delta(z^{\prime}_i,\eta_i)).
\end{eqnarray}

  Define the open sets $V_i=\Delta(z^{\prime}_{i},\eta_{i})\cap \mathring{K}_0, i\in I_0$, write $U_0\equiv V_{i_0}$, where $i_0\in I_0$ is such that  $V_{i_0}\cap \tau=\emptyset$ and let $z\in (K_0\setminus \tau)$. Since $(K_0\setminus \tau)$ is connected, we can find  a finite number of sets  $\{V_{i_k}\}_{k=0}^m$  so that $V_{i_k}\cap V_{i_{k+1}}\neq \emptyset $  and $z\in V_{i_m}$. If $m>0$, we write  $U_k=V_{i_k}\cap (\mathring{K_0}\setminus \tau), k=1,\ldots,m$; since this intersection is not necessarily connected, we choose the component of the intersection  from the condition  $U_k\cap U_{k-1}\neq \emptyset$. In particular, we have that  $c_{n}(i_k)e^{y_{1,i_{k}}(z,\epsilon_n;z_{1,i_k})}=c_{n}(i_{k+1})e^{y_{1,i_{k+1}}(z,\epsilon_n;z_{1,i_{k+1}})}, z\in U_k\cap U_{k+1}$, for all $0\leq k\leq m-1$. Hence   $(c_n(i_m)e^{y_{1,i_m}(z,\epsilon_n;z_{1,i_m})},U_m)$ is the analytic continuation of  $(c_n(i_0)e^{y_{1,i_0}(z,\epsilon_n;z_{1,i_0})},U_0)$ in the interior of the simply connected set $ K_0\setminus \tau$. If $\Delta(z^{\prime}_j,\eta_j)$ is any disk in the covering $\mathfrak{S}$ such that $\Delta(z^{\prime}_j,\eta_j)\cap U_m\neq \emptyset$ we have 
\begin{eqnarray}\label{equa3}
y_{1,i_0}(z,\epsilon_n;z_{1,i_0})=y_{1,j}(z,\epsilon_n;z_{1,j})+k_{n,j}, \quad  z\in U_m\cap \Delta(z^{\prime}_j,\eta_j), k_{n,j}\in \CC.
\end{eqnarray}
Hence, by \eqref{equa30}
\begin{eqnarray*}
\mathfrak{w}_{1,i_0}(z,\epsilon_n)=\mathfrak{w}_{1,j}(z,\epsilon_n),\quad  z\in U_m\cap \Delta(z^{\prime}_j,\eta_j).
\end{eqnarray*}
Consequently, the analytic continuation  of the function  element $(\mathfrak{w}_{1,i_0}(z,\epsilon_n),U_0)$ defines a holomorphic function $\mathfrak{w}(z,\epsilon_n)$  when   $z$ varies in the interior of the simply connected set $ K_0\setminus \tau$.
 Hence
\begin{equation}\label{equa1}
Q_{n}(z)=c_{n}e^{y(z,\epsilon_n;z_1)}, \quad \forall n>n_{0},\forall z\in  K_0\setminus \tau,
\end{equation}
where   $y(z,\epsilon_n;z_1)=\frac{1}{\epsilon_n}\int_{z_1}^z\mathfrak{w}(t,\epsilon_n)dt, \dsty n_0=\max_{i\in I_0} \{n_0(i)\}, z_1=z_{1,i_0}$, and $c_n=c_n(i_0)$.

Given $z_0\in \tau$, take $\Delta(z^{\prime}_j,\eta_{j})\subset \mathfrak{S}$   such that  $z_0\in \Delta(z^{\prime}_j,\eta_{j})$. Since $y_{1,j}(z,\epsilon_n;z_{1,j})$ is holomorphic for   $z\in \Delta(z^{\prime}_j,\eta_{j})$, using  \eqref{equa3}  we get that  for $n>n_0$,     the following limit exists and is finite 
\begin{eqnarray}\label{jumpc}
\dsty  \lim_{\begin{subarray}{c}
                           z\rightarrow z_0\\
                           z\in \tau^{+}
                          \end{subarray}}y(z,\epsilon_n;z_1)= \lim_{\begin{subarray}{c}
                           z\rightarrow z_0\\
                           z\in \tau^{+}
                          \end{subarray}}y_{1,j}(z,\epsilon_n;z_{1,j})+k_{n,j}.
\end{eqnarray}
Hence, the function $y$ is holomorphic in the interior of $K_0\setminus \tau$ and $\dsty  \lim_{\begin{subarray}{c}
                           z\rightarrow z_0\\
                           z\in \tau^{+}
                          \end{subarray}}y(z,\epsilon_n;z_1)$ exists and is finite. Therefore, by \eqref{equa1} and  \eqref{jumpc}, if we define $y(z_0,\epsilon;z_1)$  as 
\begin{eqnarray*}
y(z_0,\epsilon_n;z_1)=\lim_{\begin{subarray}{c}
                           z\rightarrow z_0\\
                           z\in \tau^{+}
                          \end{subarray}}y(z,\epsilon_n;z_1),
\end{eqnarray*}                          
where  $z_0\in \tau\cap K_0$, then $e^{-y(z,\epsilon_n;z_1)}$  is  holomorphic in the interior of $K_0$ and 
\begin{equation}\label{equa2}
Q_{n}(z)=c_{n}e^{y(z,\epsilon_n;z_1)}, \quad \forall n>n_{0},\forall z\in  K_0.
\end{equation}  
 On the other hand,  by \eqref{equa2} the function 
\begin{eqnarray}\label{frakw}
 \mathfrak{w}(z,\epsilon_n)=e^{-y(z,\epsilon_n;z_1)}Q_n^{\prime}(z)
\end{eqnarray} 
   can  be continued analytically in the interior of  $K_0$ for every $n>n_0$.  Since $Q_n$ does not have zeros outside of $\gamma_1$ for $n>n_0$ it follows that   $y$  can be continued analytically to $ext(\gamma)\setminus \tau $.  Therefore, by virtue of   \eqref{frakw},   $ \mathfrak{w}$  can also be continued analytically to $ext(\gamma)$ and we have that
\begin{equation}\label{equa5}
Q_{n}(z)=c_{n}e^{y(z,\epsilon_n;z_1)}, \quad \forall n>n_{0},\forall z\in ext(\gamma),
\end{equation} 
where  $y(z,\epsilon_n;z_1)=\frac{1}{\epsilon_n}\int_{z_1}^z\mathfrak{w}(t,\epsilon_n)dt, \mathfrak{w}(t,\epsilon_n)\in \mathcal{H}(ext(\gamma))$.

  Let $K\subset ext(\gamma)$ be an arbitrary compact subset.  Using Lemma \ref{compact} we can find  a finite covering  $\mathfrak{S}_1=\{\Delta(z^{\prime}_i,\eta_i)\}_{i\in I_1}$ by open disks 
 such that   $\dsty K \subset \bigcup_{i\in I_1}\Delta(z^{\prime}_i,\eta_i)$  and   
$$
Q_{n}(z)=c_{n}(i)e^{y_{1,i}(z,\epsilon_n;z_{1,i})}, \quad \forall n>n_{1}(i), \forall z\in  \Delta(z^{\prime}_i,\eta_i), i\in I_1,
$$
where $y_{1,i}(z,\epsilon_n;z_{1,i})=\frac{1}{\epsilon_n}\int_{z_{1,i}}^z\mathfrak{w}_{1,i}(t,\epsilon_n)dt, \mathfrak{w}_{1,i}(t,\epsilon_n)\in \mathcal{H}(\Delta(z^{\prime}_i,\eta_i))$ and 
\begin{eqnarray}\label{asympii}
\dsty
\mathfrak{w}_{1,i}(z,\epsilon_n)\sim \sum_{k=0}^{\infty}\mathfrak{b}_{1,i,k}(z)\epsilon_n^{k},
\end{eqnarray}
 as $n\rightarrow \infty$ uniformly in  $z\in \Delta(z^{\prime}_i,\eta_i)$. By \eqref{equa5}   we have that for $n>\max\{n_0,n_1(i)\}$
\begin{eqnarray}\label{asymptrestr}
\mathfrak{w}(z,\epsilon_n)=\mathfrak{w}_{1,i}(z,\epsilon_n), z\in \Delta(z^{\prime}_i,\eta_i).
\end{eqnarray}
Hence, by \eqref{asymptrestr} and using the  fact that the open covering is finite, we get
\begin{eqnarray*}
\dsty
\mathfrak{w}(z,\epsilon_n)\sim \sum_{k=0}^{\infty}\mathfrak{b}_{k}(z)\epsilon_n^{k} , \quad \mbox{for}\,  n\rightarrow\infty,  z\in K.
\end{eqnarray*}
Here  $\mathfrak{b}_{k}$ is such that 
\begin{eqnarray}\label{bfrakext}
\mathfrak{b}_{k}(z)=\mathfrak{b}_{1,i,k}(z), \quad z\in \Delta(z^{\prime}_i,\eta_i).
\end{eqnarray}

b) Follows from item b) of Proposition \ref{GL} and \eqref{bfrakext}. \lqqd

\begin{lemma}\label{preprevict}
Let $\mathfrak{w}$ be defined as in item a) of Lemma \ref{PreMainext}.  Then, $\mathfrak{w}$ and  $(\mathfrak{b}_{k})_{k\geq 0}$  are   holomorphic when $z$ varies in a neighborhood of $\infty$,   and $\dsty \lim_{z\rightarrow\infty}\mathfrak{b}_{k}(z)=0$, $\dsty \lim_{z\rightarrow\infty}\mathfrak{w}(z,\epsilon_n)=0$, for each $n>n_0$ fixed.
\end{lemma}
\proof
  By the relation 
\begin{eqnarray*}
 \mathfrak{w}(z,\epsilon_n)=\frac{\epsilon_n Q_n^{\prime}(z)}{Q_n(z)}, \quad z\in ext(\gamma), n>n_0,
\end{eqnarray*} 
we obtain that $\mathfrak{w}$ is holomorphic  in a neighborhood $V_{\infty}$ of $\infty$,  and $\dsty \lim_{z\rightarrow\infty}\mathfrak{w}(z,\epsilon_n)=0$ for fixed $n>n_0$.

We will prove by induction that $\mathfrak{b}_k$ is holomorphic when $z\in V_{\infty}$, and $\dsty \lim_{z\rightarrow\infty}\mathfrak{b}_{k}(z)=0$ for all $k\geq 0$.  Indeed,  by b) of Lemma \ref{PreMainext} one obtains that the statement is true for $k=0$.

 Assume that  $\mathfrak{b}_k$  is holomorphic  in $V_{\infty}$ for $0\leq k\leq p_0$. Since  $\dsty
\mathfrak{w}(z,\epsilon_n)\sim \sum_{k=0}^{\infty}\mathfrak{b}_{k}(z)\epsilon_n^{k} $ as $n\rightarrow \infty$ when $z$ varies in compact subsets  $K\subset ext(\gamma)$ one has for every $p\geq 0$,
$$\dsty  \mathfrak{w}(z,\epsilon_n)  =\sum_{k=0}^{p}\epsilon_n^{k}\mathfrak{b}_{k}(z)+r_{p}(z,\epsilon_n), \forall z\in ext(\gamma), $$
where $ r_{p}(z,\epsilon_n)=O(\epsilon_n^{p+1})$ as $z$ varies in $K$.  By the induction hypothesis it follows that $r_{p_0}(z,\epsilon_n)=\epsilon_{n}^{p_0+1}\mathfrak{b}_{p_0+1}(z)+r_{p_0+1}(z,\epsilon_n)$ is holomorphic in $V_{\infty}$ and $\dsty \lim_{z\rightarrow\infty}r_{p_0}(z,\epsilon_n)=0$. Hence, for all $n>n_0$,  
\begin{eqnarray}\label{2}
\mathfrak{b}_{p_0+1}(z)+\frac{r_{p_0+1}(z,\epsilon_n)}{\epsilon_{n}^{p_0+1}}=\sum_{k=1}^{\infty}\frac{a_k(\epsilon_n)}{z^k}, \quad z\in V_{\infty}, a_k(\epsilon_n)\in \CC.
\end{eqnarray}
Hence,
$$\frac{1}{2\pi\imath}\int_{|t|=R} t^{k-1}\left(\mathfrak{b}_{p_0+1}(t)+\frac{r_{p_0+1}(t,\epsilon_n)}{\epsilon_{n}^{p_0+1}}\right)dt=a_k(\epsilon_n),$$
 by taking limit as $n\rightarrow\infty$, we get  
$$ \frac{1}{2\pi\imath}\int_{|t|=R} t^{k-1}\left(\mathfrak{b}_{p_0+1}(t)+\lim_{n\rightarrow\infty}\frac{r_{p_0+1}(t,\epsilon_n)}{\epsilon_{n}^{p_0+1}}\right)dt=\lim_{n\rightarrow\infty}a_k(\epsilon_n).$$
Indeed, since 
$r_{p}(z,\epsilon_n)=o(\epsilon_n^{p}), \forall p\geq 0$, when $z$ varies in compact subsets of $\Omega$, we have that the limit in  the left--hand side exists and is finite, whence the limit in the right--hand side exists too. Therefore, 
$$ \int_{|t|=R} t^{k-1}\mathfrak{b}_{p_0+1}(t)dt=2\pi\imath a_k, \quad k\geq 1,$$
where $\dsty a_k=\lim_{n\rightarrow\infty}a_k(\epsilon_n)$. Hence by \eqref{2}, $\mathfrak{b}_{p_0+1}$ is holomorphic in $V_{\infty}$, which  completes the proof of the lemma.\lqqd

\begin{lemma}\label{existencepuisex}
For $n$  sufficiently large, we get
 \begin{itemize}
 \item[a)]  
 \begin{eqnarray}
 \begin{aligned}\label{Eqepsilons}
\epsilon_{n}^{-1}&=n+\sum_{j=0}^{\infty}  \gamma_jn^{-j},  \\
  \epsilon^k_{n}&=\sum_{j=0}^{\infty}  q_{j,k}n^{-k-j}, k>0,
\end{aligned}
\end{eqnarray}
where  $(\gamma_j)_{j=0}^{\infty}$ are defined from the following  Laurent expansion about $z=\infty$ 
 \begin{eqnarray}\label{unm1p}
\dsty f(z)=\sqrt[M]{\sum_{k=1}^{M}\rho_{k,k}(z)_k}=z+\sum_{j=0}^{\infty}\gamma_jz^{-j}
\end{eqnarray}
    and $\dsty q_{j,k}=\frac{1}{j!}\mathbb{P}^{(-k)}_j(\gamma_0,\ldots,j!\gamma_{j-1}), \gamma_0=  -\left(\frac{M-1}{2}-\frac{\rho_{M-1,M-1}}{M}\right),q_{0,k}=1$. Here  $\mathbb{P}_{n}^{(r)}$ denotes  the potential polynomials, cf. \cite[Th. B p.141] {comtet74}.

 \item[b)]  
  \begin{eqnarray*}  
n=\sum_{k=0}^{\infty}h_{k}\epsilon^{k-1}_{n}, 
\end{eqnarray*}
where  $h_0=1,h_1=-\gamma_0, h_2=-\gamma_1, \dsty h_{j+1}=-\gamma_j-\sum_{k=1}^{j-1} h_{k+1} q_{j-k,k}, j\geq 2$,  and   $(\gamma_j)_{j=0}^{\infty}$ is defined as in a).

 \end{itemize}

\end{lemma}

\proof
a)  From \eqref{unm1p} and   \eqref{epsilonn},  for $n$ sufficiently large we have 
\begin{eqnarray*}
\epsilon_{n}^{-1}&=&n-\left(\frac{M-1}{2}-\frac{\rho_{M-1,M-1}}{M}\right)+\sum_{j=1}^{\infty} \frac{\gamma_j}{n^j},  
\end{eqnarray*}
implying that
$$ \epsilon^k_{n}= \sum_{j=0}^{\infty}  q_{j,k}n^{-k-j}, k>0,$$
where $\dsty q_{j,k}=\frac{1}{j!}\mathbb{P}^{(-k)}_j(\gamma_0,\ldots,j!\gamma_{j-1}), \gamma_0= -\left(\frac{M-1}{2}-\frac{\rho_{M-1,M-1}}{M}\right), q_{0,k}=1$;  cf. \cite[Th.B p.141]{comtet74}.

b)    Assume that $n$ is large enough, then  from \eqref{epsilonn} we get
\begin{equation*} 
\sum_{m=1}^{M}\rho_{m,m}(n)_m-\frac{1}{\epsilon_n^M}=0.
\end{equation*}
By changing the variables $n=\omega$ and $u=\epsilon_n$  we obtain   the  algebraic equation, 
\begin{equation}\label{ch1}
\omega^M+\sum_{k=1}^{M-1}a_{k}\omega^{k} -\frac{1}{u^M}=0.
\end{equation}
where   $a_1,\dots ,a_{M-1}\in \CC$. Multiplying by $u^M$ and  changing of variable $w=\omega u $ in   \eqref{ch1}  one obtains 
\begin{equation*} 
F(u,w)=w^M +\sum_{k=1}^{M-1}a_{k}u^{M-k}w^{k}-1=0.
\end{equation*}

We have  that $F(0,1)=0$ and  $\dsty \frac{\partial F}{\partial w}(0,1)\neq 0$. From the implicit function  theorem  \cite[Th 3.11, Vol II]{Mark65} there exists a  neighborhood $V$ of $0$ and a unique  analytic function $\dsty w(u)=\sum_{k=0}^{\infty}u^{k}h_{k}$ such that $w(0)=1$ and  $F(u,w(u))=0, \forall u\in V$. Taking into account  that  $w=\omega u, u \neq 0$ we obtain that $\dsty \omega(u
)=\sum_{k=0}^{\infty}u^{k-1}h_{k}$ is a solution of \eqref{ch1} for    $u$ varying in a punctured 
neighborhood $V^{*}$ of $u=0$.  Hence, 
\begin{eqnarray}\label{hk}  
n=\sum_{k=0}^{\infty}h_{k}\epsilon^{k-1}_{n}, 
\end{eqnarray}
where $h_0=1$.  

Using \eqref{Eqepsilons} and \eqref{hk},
$$0=\sum_{j=0}^{\infty}  \gamma_jn^{-j}+h_1+\sum_{k=1}^{\infty}h_{k+1}\sum_{j=0}^{\infty}  q_{j,k}n^{-k-j},$$
after some straightforward transformations we obtain 
$$0=\gamma_0+h_1+\sum_{j=1}^{\infty}\left(\gamma_j+\sum_{k=1}^{j} h_{k+1} q_{j-k,k}\right)n^{-j},$$
yielding  $h_1=-\gamma_0, h_2=-\gamma_1$, and  $\dsty h_{j+1}=-\gamma_j-\sum_{k=1}^{j-1} h_{k+1} q_{j-k,k}, j\geq 2$. 
   \lqqd

\begin{lemma}\label{previct}
Let    $\tau$ be  any Jordan arc  from  $-\infty$ to $\omega\in \supp \mu$ such that $]-\infty,p]\subset \tau$ for some $p\in \RR$  and let   $\gamma$ be a closed Jordan curve enclosing  $ \supp\mu$. Then,  for $n$ sufficiently large, we get 
$$
Q_n(z)=e^{ \frac{1}{\epsilon_n}y(z,\epsilon_n)}, \forall z\in ext(\gamma), 
$$
where $\dsty
y(z,\epsilon_n)\sim \sum_{k=0}^{\infty}\Phi_{k}(z)\epsilon_n^{k}, $ as $n\rightarrow \infty$ when $z$ varies in compact subsets  $K\subset ext(\gamma)$ and $\Phi_k(z)$ is the primitive of the function $\mathfrak{b}_k$ in $ext(\gamma)\setminus \tau$   such that $\dsty \lim_{z\rightarrow\infty}\Phi_k(z)-h_k\ln z=0$, $\dsty \Phi_k(z)=\lim_{\begin{subarray}{c}
                           z\rightarrow z_0\\
                           z\in \tau^{+}
                          \end{subarray}}\Phi_k(z)$,  when  $z_0\in \tau\cap ext(\gamma)$, being $h_k$ as in item  b) of Lemma \ref{existencepuisex}.
\end{lemma}

\proof
  By item a) of Lemma \ref{PreMainext}  there exist  $z_1\in ext(\gamma), n_0, y$, and a constant $c_n$  such that    
$$Q_{n}(z)=c_{n}e^{\dsty y(z,\epsilon_n;z_1)}, \forall n>n_{0}, \quad \forall z\in ext(\gamma),$$
where  $y(z,\epsilon_n;z_1)=\frac{1}{\epsilon_n}\int_{z_1}^z\mathfrak{w}(t,\epsilon_n)dt, \mathfrak{w}(t,\epsilon_n)\in \mathcal{H}(ext(\gamma))$,   $\dsty
\mathfrak{w}(z,\epsilon_n)=\sum_{k=0}^{p}\mathfrak{b}_{k}(z)\epsilon_n^{k}+\mathfrak{r}_p(z,\epsilon_n)\epsilon_n^{p+1}$, and  as $n\rightarrow\infty$,
\begin{eqnarray}\label{remcond}
\mathfrak{r}_{p}(z,\epsilon_n)=O(1),
\end{eqnarray}
 for   $z$ varying in compact subsets  $K\subset ext(\gamma)$, and  any fixed $ p\geq 0$.

 Therefore,  for $n$  fixed,
$$
\dsty \lim_{ z\rightarrow \infty} \frac{e^{ y(z,\epsilon_n;z_1)}}{z^n}=\frac{1}{c_n}.
$$
Hence, for every $p\geq 0$ fixed
\begin{equation} \label{limit}
\dsty \lim_{ z\rightarrow \infty}\exp\left(\dsty
\sum_{k=0}^{p}\epsilon_n^{k-1}\int_{z_1}^z\mathfrak{b}_k(t)dt+\epsilon_n^{p}\int_{z_1}^z\mathfrak{r}_{p}(t,\epsilon_n)dt-n\ln z\right)=\frac{1}{c_n}.
\end{equation}

From item b) of Lemma \ref{existencepuisex},  every sufficiently large   $n\in \NN$  can be uniquely  expressed  as
\begin{eqnarray}\label{n}
n=\sum_{k=0}^{p}h_{k}\epsilon_n^{k-1}+\epsilon_n^{p}g_p(\epsilon_n), \forall p\in \NN\cup\{0\},
\end{eqnarray}
where   $h_k$ is defined as in item b) of Lemma \ref{existencepuisex} and 
\begin{eqnarray}\label{gp}
g_p(\epsilon_n)=O(1), \quad \mbox{as} \quad  n\rightarrow \infty, 
\end{eqnarray}
for all $p\geq 0$ fixed. 

From   \eqref{limit}  and \eqref{n}, for all $p\geq 0$ fixed, we get 
\begin{eqnarray} \label{limita}
\dsty \lim_{ z\rightarrow \infty}\exp\left(\dsty \epsilon_{n}^{-1}\left(\int_{z_1}^z\mathfrak{b}_0(t)dt -\ln z\right)+
\sum_{k=1}^{p}\left(\int_{z_1}^z\mathfrak{b}_k(t)dt
-h_{k}\ln z\right)\epsilon_{n}^{k-1}+
\left(\int_{z_1}^z\mathfrak{r}_p(t,\epsilon_n)dt-g_p(\epsilon_n)\ln z\right)\epsilon_n^{p}\right)=\frac{1}{c_n}.
\end{eqnarray}

%Consider that  $\Phi_k(z)=\int_{z_1}^z\mathfrak{b}_k(t)dt+d_k$ and 

Let $h_k, k\geq 0$ be defined as in item b) of Lemma \ref{existencepuisex}. We want to prove that for all $p\geq 0$, there exists a primitive of the function $\mathfrak{b}_p$  in $ext(\gamma)\setminus \tau$ such that   $\dsty \lim_{ z\rightarrow \infty}\Phi_p(z)-h_p\ln z=0$ and  
\begin{eqnarray}\label{wewant}
Q_n(z)=\exp\left(\dsty  
\sum_{k=0}^{p}\Phi_k(z)\epsilon_{n}^{k-1}+\Psi_p(z,\epsilon_n)\right),
\end{eqnarray}
where $\Psi_p(z,\epsilon_n)=O(\epsilon^p_n)$ as $n\rightarrow \infty$ when $z$ varies in compact subsets  $K\subset ext(\gamma)$.

The argument follows  by induction in $p$. For   $p=0$, by  \eqref{limit}  and \eqref{n}, we get 
\begin{eqnarray} \label{limitini}
\dsty \lim_{ z\rightarrow \infty}\exp\left(\dsty \epsilon_{n}^{-1}\left(\int_{z_1}^z\mathfrak{b}_0(t)dt -\ln z\right)+\int_{z_1}^z\mathfrak{r}_0(t,\epsilon_n)dt-g_0(\epsilon_n)\ln z\right)=\frac{1}{c_n}.
\end{eqnarray}

By item b) of Lemma \ref{PreMainext} we have that   $\dsty \mathfrak{b}_0(z)=w_1(z), z\in ext(\gamma)$. Let  $\Phi_0$ be  the primitive of $ \mathfrak{b}_0$ which asymptotically  coincides with $\ln z$ near of $\infty$. We have that 
\begin{eqnarray}\label{primi0}
\begin{aligned}
\Phi_0(z)&=\int_{z_1}^z\mathfrak{b}_0(t)dt+\Phi_0(z_1), \quad  z\in ext(\gamma)\setminus \tau,\\
\Phi_0(z_0)&=\lim_{\begin{subarray}{c}
                           z\rightarrow z_0\\
                           z\in \tau^{+}
                          \end{subarray}}\Phi_0(z),\quad   z_0\in \tau\cap ext(\gamma),
\end{aligned}
\end{eqnarray}

On the other hand, the  relation \eqref{limitini} gives 
$$\dsty \lim_{  z\rightarrow \infty} \int_{z_1}^z\mathfrak{r}_0(t,\epsilon_n)dt-g_0(\epsilon_n)\ln z=d_n\in \CC,$$
therefore, by taking 
\begin{eqnarray}\label{primi01}
\begin{aligned}
\Psi_0(z,\epsilon_n)&=\int_{z_1}^z\mathfrak{r}_0(t,\epsilon_n)dt-d_n,\quad  z\in ext(\gamma)\setminus \tau,\\
\Psi_0(z_0,\epsilon_n)&=\lim_{\begin{subarray}{c}
                           z\rightarrow z_0\\
                           z\in \tau^{+}
                          \end{subarray}}\Psi_0(z,\epsilon_n),\quad   z_0\in \tau\cap ext(\gamma),
\end{aligned}
\end{eqnarray}
and by using  \eqref{primi0} and \eqref{primi01} we obtain 
$$
\dsty \lim_{ z\rightarrow \infty}\exp\left(\dsty \epsilon_{n}^{-1}\left(\Phi_0(z)-\ln z\right)+\Psi_0(z,\epsilon_n)-g_0(\epsilon_n)\ln z\right)=1.
 $$
Notice that by virtue of \eqref{remcond} and \eqref{gp}, we obtain 
$$\Psi_0(z,\epsilon_n)=O(1),$$
as $n\rightarrow\infty$ and for  $z$ varying in compact subsets of $ext(\gamma)$. Hence, relation \eqref{wewant}  holds for $p=0$.

Assume that \eqref{wewant} holds for $p=p_0$. From the relation 
$$Q_n=c_n\exp\left(\dsty \sum_{k=0}^{p_0}\Phi_{k}(z) \epsilon_{n}^{k-1}+\epsilon_{n}^{p_0}\int_{z_1}^z\mathfrak{b}_{p_0+1}(t)dt+\epsilon_{n}^{p_0+1}\int_{z_1}^z\mathfrak{r}_{p_0+1}(t,\epsilon_n)dt\right), \forall z\in ext(\gamma), n>n_0, $$
and from the fact that  
 \begin{multline}\label{p0+1}
\lim_{z\rightarrow \infty }\exp\left(\dsty 
\sum_{k=0}^{p_0}\left(\Phi_k(z)
-h_{k}\ln z\right)\epsilon_{n}^{k-1}+\left(\int_{z_1}^z\mathfrak{b}_{p_0+1}(t)dt
-h_{p_0+1}\ln z\right)\epsilon_{n}^{p_0}+\right.\\
\left.
\left(\int_{z_1}^z\mathfrak{r}_{p_0+1}(t,\epsilon_n)dt-g_{p_0+1}(\epsilon_n)\ln z)\right)\epsilon_{n}^{p_0+1}\right)=
\frac{1}{c_n}, 
\end{multline}
it follows that 
 \begin{equation}\label{indtes}
\lim_{\begin{subarray}{c}
     z\rightarrow \infty\\
   z\in  ext(\gamma)\setminus \tau
                          \end{subarray}}\left(\int_{z_1}^z\mathfrak{b}_{p_0+1}(t)dt
-h_{p_0+1}\ln z\right)\epsilon_{n}^{p_0}+
\left(\int_{z_1}^z\mathfrak{r}_{p_0+1}(t,\epsilon_n)dt-g_{p_0+1}(\epsilon_n)\ln z)\right)\epsilon_{n}^{p_0+1}=d_n^{\prime}\in \CC.
\end{equation}

By Lemma \ref{preprevict}, there exists   a neighborhood  $V_{\infty}$ of $\infty$ such that for $z\in  V_{\infty}\setminus \tau$  
\begin{eqnarray}
\label{lau1}\int_{z_1}^z\mathfrak{b}_{p_0+1}(t)dt&=&\sum_{k=2}^{\infty}\frac{\alpha_{-k,p_0+1}}{z^{k-1}}+\alpha_{-1,p_0+1}\ln z+\alpha_{0,p_0+1} ,\\
\label{lau2}\int_{z_1}^z\mathfrak{r}_{p_0+1}(t,\epsilon_n)dt&=&\sum_{k=2}^{\infty}\frac{m_{-k,p_0+1}(\epsilon_n)}{z^{k-1}}+m_{-1,p_0+1}(\epsilon_n)\ln z +m_{0,p_0+1}(\epsilon_n).
\end{eqnarray}

Hence, by  \eqref{indtes}, \eqref{lau1} and \eqref{lau2}  we deduce that 
\begin{multline}\label{1lim}
\lim_{\begin{subarray}{c}
                           z\rightarrow \infty\\
                           z\in  V_{\infty}\setminus \tau
                          \end{subarray}}\sum_{k=2}^{\infty}\frac{\alpha_{-k,p_0+1}\epsilon_{n}^{p_0}+m_{-k,p_0+1}(\epsilon_n)\epsilon_{n}^{p_0+1}}{z^{k-1}}+((\alpha_{-1,p_0+1}-h_{p_0+1})\epsilon_{n}^{p_0}+
                   (m_{-1,p_0+1}(\epsilon_n)-g_{p_0+1}(\epsilon_n))\epsilon_{n}^{p_0+1})\ln z+\\
                   (\alpha_{0,p_0+1}\epsilon^{p_0}_n +m_{0,p_0+1}(\epsilon_n) \epsilon^{p_0+1}_n),
\end{multline}
for all fixed $n>n_0$   exists and is finite.

Since $\mathfrak{r}_{p_0+1}$ is holomorphic for $z\in V_{\infty}$ one has 
$\dsty  m_{-1,p_0+1}(\epsilon_n)2\pi\imath=\int_{|t|=R} \mathfrak{r}_{p_0+1}(t,\epsilon_n)dt$, being $R$ large enough. 
Hence by \eqref{remcond},
\begin{eqnarray}\label{2lim}
 \lim_{n\rightarrow\infty}m_{-1,p_0+1}(\epsilon_n),
\end{eqnarray}
 exists and is finite.

On the other hand, by \eqref{1lim}
\begin{eqnarray}\label{precoeffprimitiver}
 \alpha_{-1,p_0+1}-h_{p_0+1}=-(m_{-1,p_0+1}(\epsilon_n)-g_{p_0+1}(\epsilon_n))\epsilon_{n},
\end{eqnarray}
therefore, by  \eqref{gp} and \eqref{2lim} 
$$ \lim_{n\rightarrow\infty}-(m_{-1,p_0+1}(\epsilon_n)-g_{p_0+1}(\epsilon_n))\epsilon_{n}=0=\alpha_{-1,p_0+1}-h_{p_0+1},$$
that is 
\begin{eqnarray}
\label{coeffprimitiveb} \alpha_{-1,p_0+1}&=& h_{p_0+1}.
\end{eqnarray} 
Hence, from \eqref{precoeffprimitiver}
\begin{eqnarray}
\label{coeffprimitiver} m_{-1,p_0+1}(\epsilon_n)&=& g_{p_0+1}(\epsilon_n).
\end{eqnarray} 

By \eqref{coeffprimitiveb} and \eqref{coeffprimitiver} we deduce that there exist a primitive  $\Phi_{p_0+1}$  of the function $\mathfrak{b}_{p_0+1}$ and a primitive $\Psi_{p_0+1}$ of the function $\mathfrak{r}_{p_0+1}$ in $ext(\gamma)\setminus \tau $ such that 
\begin{eqnarray*}
\nonumber \dsty \lim_{  z\rightarrow \infty }\Phi_{p_0+1}(z)-h_{p_0+1}\ln z&=&0,\\
\nonumber \dsty \lim_{  z\rightarrow \infty }  \Psi_{p_0+1}(z,\epsilon_n)-g_{p_0+1}(\epsilon_n)\ln z&=&0,
\end{eqnarray*} 
hence, if $\dsty \Phi_{p_0+1}(z_0)=\lim_{\begin{subarray}{c}
                           z\rightarrow z_0\\
                           z\in \tau^{+}
                          \end{subarray}}\Phi_{p_0+1}(z)$ and $\dsty  \Psi_{p_0+1}(z_0,\epsilon_n)=\lim_{\begin{subarray}{c}
                           z\rightarrow z_0\\
                           z\in \tau^{+}
                          \end{subarray}}\Psi_{p_0+1}(z,\epsilon_n)$ when  $z_0\in \tau\cap ext(\gamma)$ we obtain   for all fixed $n>n_0$, 

$$
\lim_{z\rightarrow \infty }e^{\dsty 
\sum_{k=0}^{p_0}\left(\Phi_k(z)
-h_{k}\ln z\right)\epsilon_{n}^{k-1}+\left(\Phi_{p_0+1}(z)
-h_{p_0+1}\ln z\right)\epsilon_{n}^{p_0}+
\left( \Psi_{p_0+1}-g_{p_0+1}(\epsilon_n)\ln z)\right)\epsilon_{n}^{p_0+1}}=1,
$$
notice that by \eqref{remcond} and \eqref{gp}
$$\Psi_{p_0+1}(z,\epsilon_n)=O(1),$$
as $n\rightarrow\infty$ when  $z$ varies in compact subsets of $ext(\gamma)$. Hence 
\eqref{wewant} holds for $p=p_0+1$. This completes the proof of the  lemma.
 \lqqd

\proof{(\emph{of Theorem \ref{Main}})}

Notice that  by \cite[Th.3]{BerRull02}, the set $\supp\mu$ is a finite  tree.  It suffices to consider any Jordan curve $\gamma$ enclosing  $\supp \mu$ and show  the asymptotic relation in  compact subsets of $ext(\gamma)$.
 By Lemma \ref{previct}, 
\begin{eqnarray}\label{Qnpre}
Q_n(z)=\exp\left(\dsty \sum_{k=0}^{p}\Phi_{k}(z) \epsilon_{n}^{k-1}+\Psi_{p}(z,\epsilon_n)\epsilon_n^{p}\right), \forall z\in ext(\gamma), \forall p\in \NN\cup\{0\},
\end{eqnarray}
where   $\Phi_k$ is the primitive of the function $\mathfrak{b}_k$  in $ext(\gamma)\setminus \tau$     such that $\dsty \lim_{z\rightarrow\infty}\Phi_k(z)-h_k\ln z=0$, $h_{k}$ as in   b) of  Lemma \ref{existencepuisex}, and $\Psi_{p}(z,\epsilon_n)=O(1)$ as $n\rightarrow\infty$, for   $z$ varying in compact subsets  $K\subset ext(\gamma)$,  $\forall p\geq 0$ fixed.

Using     a) of  Lemma \ref{existencepuisex}, write \eqref{Qnpre} for $p>1$  as 
\begin{multline*}
  Q_n(z)=\\
  \exp\left(\dsty \left(n-\left(\frac{M-1}{2}-\frac{\rho_{M-1,M-1}}{M}\right)\right)\Phi_{0}(z)+\Phi_{1}(z)+\Phi_0(z)\sum_{j=1}^{p-1} \frac{\gamma_j}{n^j}+ \sum_{k=1}^{p-1}\sum_{j=k}^{p-1}\frac{q_{j-k,k}}{n^j}\Phi_{k+1}(z) +\Psi^{*}_{p}(z,n)\frac{1}{n^{p}}\right)=\\
  \exp\left(\dsty \dsty n\Phi_{0}(z)-\left(\frac{M-1}{2}-\frac{\rho_{M-1,M-1}}{M}\right)\Phi_{0}(z)+\Phi_{1}(z)+\sum_{j=1}^{p-1} \Phi^*_j(z)\frac{1}{n^j}+\Psi^*_{p}(z,n)\frac{1}{n^{p}}\right), z\in ext(\gamma),
\end{multline*}
where $\dsty \Phi^*_j(z)=\gamma_j\Phi_0(z)+ \sum_{k=1}^{j}q_{j-k,k}\Phi_{k+1}(z)$. Since $\Psi^*_{p}(z,\epsilon_n)=O(1)$ as $n\rightarrow\infty$, for   $z$ varying in compact subsets  $K\subset ext(\gamma)$, we can write 
\begin{eqnarray}\label{almost}
Q_n(z)=e^{\dsty \dsty n\Phi_{0}(z)-\left(\frac{M-1}{2}-\frac{\rho_{M-1,M-1}}{M}\right)\Phi_{0}(z)+\Phi_{1}(z)}\left(1+\sum_{j=1}^{p-1}\frac{C_j(z)}{n^j}+O(1/n^{p})\right),
\end{eqnarray}
where $\dsty C_j(z)=\frac{1}{j!}{\bf Y}_j(\Phi^*_1(z),2\Phi^*_2(z),\ldots,j!\Phi^*_j(z))$.

To prove  that $C_j\in\mathcal{H}(\Omega), \forall j\in \NN$,  we first  prove  that $\exp(\dsty \Phi_0)$ and  $\exp\left(\dsty \left(\frac{M-1}{2}-\frac{\rho_{M-1,M-1}}{M}\right)\Phi_{0}+\Phi_{1}\right)$ belong to $\mathcal{H}(\Omega)$. We may assume that $\Phi_0(z)\neq \ln z$, or equivalently that $\rho_m(z)\neq z^M$, otherwise there is nothing to prove.  By  Lemma \ref{preprevict},  $\mathfrak{b}_0\in \mathcal{H}(V_{\infty})$    and $\dsty \lim_{z\rightarrow\infty}\mathfrak{b}_{k}(z)=0$, being $V_{\infty}$ a neighborhood of $\infty$. Since  $\Phi_0$ is the primitive of the function $\mathfrak{b}_0$ in $ext(\gamma)$   such that $\dsty \lim_{z\rightarrow\infty}\Phi_0(z)-\ln z=0$  we have that   $e^{\dsty \Phi_0}$ admits a Laurent  expansion
$$e^{\dsty \Phi_0(z)}=\beta z+\sum_{k=0}^{\infty}\frac{\beta_k}{z^k}, \quad z\in V_{\infty},$$
 write
\begin{eqnarray}\label{mono}
e^{\dsty \Phi_0(z)}-\beta z=\sum_{k=0}^{\infty}\frac{\beta_k}{z^k}, \quad z\in V_{\infty}.
\end{eqnarray}
Take  $\widetilde{\Omega}=\mathcal{T}(\Omega)$, where $\mathcal{T}(z)=\frac{1}{z}$ and  define the boundary point $\widetilde{\omega}$ of $\widetilde{\Omega}$  by the condition $\dsty |\widetilde{\omega}|=\min_{\omega\in \partial \widetilde{\Omega}}|\omega| $. Taking  a cut in $\widetilde{\Omega}$ using the  line segment $\widetilde{\mathfrak{l}}$ from $0$ to $\widetilde{\omega}$, we see  that $\mathfrak{l}=\mathcal{T}^{-1}(\widetilde{\mathfrak{l}})$ is a Jordan arc in the region $\Omega$. By \eqref{mono}, we have that  the left--hand side  is holomorphic in  $V_{\infty}$, hence  by  \cite[Th. 8.5, p.269 \& Cor. p.272, Vol. III]{Mark65} it can be prolonged analytically in  the simply connected open set  $V_{\infty}\bigcup(\Omega\setminus \mathfrak{l})$. Now, by considering $z\mapsto \frac{1}{z}$ we see that the right--hand side of \eqref{mono} is a convergent power series in the disk $\Delta(0,|\widetilde{\omega}|)\subset \widetilde{\Omega}$. Therefore $e^{\dsty \Phi_0(z)}\in \mathcal{H}(\Omega)$.

On the other hand, by \eqref{Qnpre} one has 
\begin{eqnarray*}
\dsty \lim_{n\rightarrow \infty}\frac{Q_n(z)}{e^{\dsty n\Phi_0(z)}}=\exp\left(\dsty \left(\frac{M-1}{2}-\frac{\rho_{M-1,M-1}}{M}\right)\Phi_{0}(z)+\Phi_{1}(z)\right), 
\end{eqnarray*}
in compact subsets of $\Omega$. Since $e^{\dsty \Phi_0}\in \mathcal{H}(\Omega)$,  it follows from the Weierstrass theorem on uniformly convergent sequences of holomorphic functions cf. \cite[Th.15.8  p.330 Vol I]{Mark65} that $\dsty \exp\left(\dsty \left(\frac{M-1}{2}-\frac{\rho_{M-1,M-1}}{M}\right)\Phi_{0}+\Phi_{1}\right) \in \mathcal{H}(\Omega)$. 

Now, from the fact that $\exp\left(\dsty \Phi_0\right), \exp\left(\dsty \left(\frac{M-1}{2}-\frac{\rho_{M-1,M-1}}{M}\right)\Phi_{0}+\Phi_{1}\right) \in \mathcal{H}(\Omega)$ and from \eqref{almost} we have that in compact subsets of $\Omega$
$$C_1(z)=\lim_{n\rightarrow\infty}\frac{nQ_n(z)}{\exp\left(\dsty \dsty n\Phi_{0}(z)-\left(\frac{M-1}{2}-\frac{\rho_{M-1,M-1}}{M}\right)\Phi_{0}(z)+\Phi_{1}(z)\right)}-n-O\left(\frac{1}{n}\right),$$
hence, $C_1\in \mathcal{H}(\Omega)$. A straightforward argument by induction  shows  that $C_j\in\mathcal{H}(\Omega), \forall j\in \NN$, which  completes the proof of the theorem. \lqqd

%%%%%%%%%%%%%%%%%%%%%%%%%%%%%%%%%%%%%%%%%%%%%%%%%%%%%%%%%%%%%%%%%%%%%%%%%%%%%%%%%%%%%%%%%%%%%%%%%%%%%%%%%%%%%%%%%%%%%%%%%%%%%%%%%%%%%%%%%%%%%%

\section{Stokes phenomenon}\label{Stk}

Recall that Stokes phenomenon or Stokes discontinuity  refers to the discontinuous change in an  asymptotic expansion through certain rays. It was  named after George Stokes  who first discovered it \cite{Stk64,Stbook}. Stokes phenomenon   has been a topic of interest for over a century e.g.  \cite{Ablowitz03,Ding, PaW95}.

 From the results about the asymptotic expansion, we see that for $\rho_M(z)=z^M$ the asymptotic expansion reads as
$$\dsty Q_n(z)\sim z^ne^{\dsty \frac{1}{M}\sum_{k=1}^{M-1}\frac{\rho_{M-1,M-1-k}}{kz^k}}$$
in compact subset of $\CC\setminus \{0\}$. Hence, if $\rho_{M-1}(z)=z^{M-1}$ or $\rho_{M-1}(z)=0$ then $\dsty Q_n(z)\sim z^n$ and we have that $\infty$ is the unique singularity, therefore we do not have the  Stokes phenomenon in $\CC$. If $\rho_{M-1}(z)\neq z^{M-1}$ we have a discontinuity in the asymptotic expansion at $z=0$ and $z=\infty$.  On the other hand, when  $\rho_M(z)\neq z^M$  we have that the Stokes phenomenon  occurs at the tree defined by  $\supp \mu$. The description of the Stokes phenomena of the present asymptotic expansion by means of the   Stokes curves, cf. \cite{Fed93,H15vir,MT17} is certainly  an interesting and non-trivial  problem  worth to be studied further.

\section{Applications}\label{Ap}

\subsection{A conjecture on a class of differential operators}
Let $P_M$ be a monic  polynomial of degree $M$ and consider the exactly solvable $M$th--order linear differential operator
\begin{eqnarray}\label{MS}
\LL^{(M)}[v]=\frac{d^M}{dz^M}\left( P_Mv \right). 
\end{eqnarray}

This operator; has been considered  in  \cite{masshap01},  the authors formulate a series of conjectures about  polynomial eigenfunctions $Q_n$, which we will  assume monic without loss of generality. Most of these conjectures have been settled in \cite{Berg07,BerRull02,BeRuSh04}. In particular, we are interested in  Conjecture $11$ of loc. cit.
\begin{conjecture}[\cite{masshap01}]
Let $BC$ be  a branch cut in $\CC$ consisting of rays from all distinct roots of $P_M$ to $\infty$ so that they do not intersect each other or the tree defined by $\supp \mu$. Select the branch of $\frac{1}{\sqrt[M]{P_M(z)}}$ in $\Omega_{P_M}=\CC\setminus BC$ which asymptotically coincides with $\frac{1}{z}$ near $\infty$ and choose $z_0\in \Omega_{P_M}$. Then, 
$$\dsty \lim_{n\rightarrow\infty}\frac{Q_n(z)}{e^{n\Psi_{P_M}(z)}}= \left( \frac{e^{\Psi_{P_M}(z)}}{\sqrt[M]{P_M(z)}} \right)^{\frac{M+1}{2}},$$
uniformly in compact subsets of  $\CC\setminus \supp \mu$, where   $\dsty \Psi_{P_M}(z)=\int_{z_0}^z \frac{1}{\sqrt[M]{P_M(t)}}dt$. 
\end{conjecture}

It should be remarked  that in  \cite{masshap01} the authors claim that this conjecture is proven in \cite{BerRull02}, however in  this publication the authors do not study  the strong asymptotic behavior of the sequence of the  polynomial eigenfunctions. 

The  result in Theorem \ref{Main} can be applied to the operator defined by \eqref{MS}. Indeed, let $\Phi_0$ be as in Theorem \ref{Main}, a straightforward calculation shows  that 
\begin{eqnarray*}
\Phi_1(z)=-\frac{M+1}{2M}\ln P_M(z), 
\end{eqnarray*}
hence,
$$ \lim_{n\rightarrow\infty}\frac{Q_n(z)}{e^{n\Phi_0(z)}}= \left( \frac{e^{\Phi_0(z)}}{\sqrt[M]{P_M(z)}} \right)^{\frac{M+1}{2}},$$
in compact subsets of $\CC\setminus \supp \mu$.  Observe that the results do not coincide in the factor $\sqrt[M]{P_M(z)}$ since the branch cuts are different.

\subsection{Asymptotic behavior of the monic eigenpolynomials of a fourth order differential operator}

As a concrete example, consider the  linear differential operator of order four, \cite{JuKwLee97}
$$\LL^{(M)}[v]=(z^2-1)^2u^{(4)}+4z(z^2-1)v^{(3)}+2(z-1)((1+2\,c)z+2\,c+3)v^{(2)},$$
denote by $(Q_n)_{n=0}^{\infty}$   the sequence of its  monic polynomial eigenfunctions with eigenvalues $\lambda_n=n(n-1)(n^2-n+4\,c)$. It is known that the sequence  is orthogonal with respect to the Sobolev inner product
 $$\langle P,Q\rangle=P(1)\overline{Q}(1)+ \frac{1}{c}\,P^{\prime}(1)\overline{Q}^{\prime}(1)+\int_{-1}^{1}P^{\prime}\overline{Q}^{\prime}dx,\quad c >0.$$

The results of the preceding section can be applied to study the strong asymptotic behavior of the sequence $(Q_n)_{n=0}^{\infty}$.   Let  $\gamma$ be  any  closed Jordan curve enclosing $[-1,1]$. From Theorem \ref{Main} we obtain immediately 
$$\Phi_0(z)=\int^z \frac{d\zeta}{\sqrt{\zeta^2-1}}=\ln\left(\frac{\varphi(z)}{2}\right), \quad z\in ext(\gamma),$$
where $\varphi(z)=z+\sqrt{z^2-1}$, here  we take the branch of  $\varphi$ for which $|\varphi(z)|>1$ whenever $z \in \CC \setminus [-1,1]$. In a similar way,
\begin{eqnarray*}
\nonumber \Phi_1(z)&=&\int^z  \left(\frac{3\rho_{M}^{\prime}(\zeta)}{8(\zeta^2-1)^2}-\frac{z}{\zeta^2-1}\right)d\zeta \\
\nonumber &=&\frac{1}{4}  \ln\left(z^2-1\right), \quad z\in  ext(\gamma).
\end{eqnarray*}
Therefore,
\begin{itemize}
\item [a)] $\dsty  Q_n(z)\sim\sqrt{2}\left(\frac{\varphi(z)}{2}\right)^{n}\frac{\sqrt[4]{z^2-1}}{\sqrt{\phi(z)}}\left(1+\frac{C_1(z)}{n}+\frac{C_2(z)}{n^2}+\ldots\right),$
\item [b)]$\dsty\lim_{n\rightarrow \infty} \frac{Q_{n+1}(z)}{Q_{n}(z)}= \frac{\varphi(z)}{2}$,
\item [c)] $\dsty\lim_{n\rightarrow \infty} \sqrt[n]{|Q_n(z)|}= \frac{|\varphi(z)|}{2},$
\end{itemize}
uniformly in compact subsets  $K\subset ext(\gamma)$.

\section{Concluding Remarks}

In the present manuscript, we obtain an asymptotic expansion in the sense of Poincar\'e for  the  polynomial eigenfunctions (eigenpolynomials) of  exactly solvable differential operators. The result is derived basically from a classical theorem due to Sibuya on the existence of asymptotic solutions   for nonlinear systems of first order differential equations depending on a small parameter and from another important publication  \cite{BerRull02}. The study of asymptotic expansions has drawn a great deal of attention for over a century, we would like to pose the following problem, asked by an anonymous reviewer,  to be considered further,

\begin{problem}

Is it possible to extend the results/methods of the paper to solutions of exactly solvable operators other than eigenpolynomials?
\end{problem}

Worth to say that some key ingredients such as the ones given in Lemma \ref{SuperSt} probably would not be valid in a general case. The study of the Stokes geometry for the analysis of the Borel summability of  WKB solutions seems a promising path to attack this problem.

%%%%%%%%%%%%%%%%%%%%%%%%%%%%%%%%%%%%%%%%%%%%%%%%%%%%%%%%%%%%%%%%%%%%%%%%%%%%%%%%%%%%%%%%%%%%%%%%%%%%%%%%%%%%%%%%%%%%%%%%%%%%%%%%%%%%%%%%%%%%%%%%%%%%%%5555

 \section{Acknowledgments}

The author  acknowledges    to  Boris Shapiro, Rikard B\"ogvad, and the anonymous reviewers for the careful reading and  helpful comments which improved the quality of the manuscript.

\bibliography{J_Borrego}

\end{document}